\documentclass[reqno,a4paper,10pt,twoside]{amsart}
\usepackage{amsmath,amssymb,epsf,epic}
\usepackage{times}
\usepackage{mathrsfs}
\usepackage{latexsym,amsopn}
\usepackage{amsfonts,amsbsy,amscd,stmaryrd}
\usepackage{paralist}
\usepackage{geometry}
\numberwithin{equation}{section}

\renewcommand{\ker}{{\rm Ker}}
\newcommand{\bra}{\langle} 
\newcommand{\ket}{\rangle}

\renewcommand\Im{{\rm Im\,}}
\renewcommand\Re{{\rm Re\,}}
\renewcommand{\min}{{\rm min}}
\renewcommand{\max}{{\rm max}}

\newcommand{\cl}{{\mathcal L}}

\newcommand{\e}{\mathrm{e}}
\renewcommand{\c}{\mathrm{c}}

\renewcommand{\i}{{\rm i}}
\renewcommand{\d}{{\rm d}}

\renewcommand{\sp}{{\rm sp\,}}

\newcommand{\rs}{{\rm rs\,}}
\newcommand{\aam}{A_\alpha^{\min}}
\newcommand{\aaM}{A_\alpha^{\max}}
\newcommand{\m}{{\min}}
\newcommand{\M}{{\max}}

\newcommand{\grintl}{[\kern-.18em [}
\newcommand{\grintr}{]\kern-.18em ]}
\newcommand{\ds}{\displaystyle}
\def\slim{{\rm s-}\lim}

\def\Num{{\rm Num\,}}
\def\Ran{{\rm Ran\,}}
\def\ccf{C_\c^\infty}

\newcounter{resultcounter}[section]
\renewcommand{\theresultcounter}{\arabic{section}.\arabic{resultcounter}}

\newtheorem{theorem}[resultcounter]{Theorem}
\newtheorem{lemma}[resultcounter]{Lemma}
\newtheorem{proposition}[resultcounter]{Proposition}
\newtheorem{corollary}[resultcounter]{Corollary}
\newtheorem{definition}[resultcounter]{Definition}
\newtheorem{remark}[resultcounter]{Remark}

\newtheorem{example}[resultcounter]{Example}

\newcommand\bep{\begin{proposition}}
\newcommand\eep{\end{proposition}}
\newcommand\ber{\begin{remark}}
\newcommand\eer{\end{remark}}
\newcommand\bel{\begin{lemma}}
\newcommand\eel{\end{lemma}}
\newcommand\bet{\begin{theorem}}
\newcommand\eet{\end{theorem}}
\newcommand\bex{\begin{example}}
\newcommand\eex{\end{example}}
\newcommand\bed{\begin{definition}}
\newcommand\eed{\end{definition}}
\newcommand\bea{\begin{assumption}}
\newcommand\eea{\end{assumption}}
\newcommand\bec{\begin{corollary}}
\newcommand\eec{\end{corollary}}
\newcommand{\beq}{\begin{equation}}
\newcommand{\eeq}{\end{equation}}
\newcommand{\beqa}{\begin{eqnarray}}
\newcommand{\eeqa}{\end{eqnarray}}

\def\one{{\mathchoice {\rm 1\mskip-4mu l} {\rm 1\mskip-4mu l} {\rm
      1\mskip-4.5mu l} {\rm 1\mskip-5mu l}}} 

 \def\cB{{\mathcal B}} 
\def\cD{{\mathcal D}} \def\cE{{\mathcal E}} \def\cF{{\mathcal F}}
\def\cG{{\mathcal G}} \def\cH{{\mathcal H}} 
\def\cJ{{\mathcal J}} \def\cK{{\mathcal K}} \def\cL{{\mathcal L}}
\def\cM{{\mathcal M}}  
 \def\cQ{{\mathcal Q}}

\newcommand{\R}{{\mathbb R}}
\newcommand{\N}{{\mathbb N}}

\newcommand{\C}{{\mathbb C}}
\newcommand{\Z}{{\mathbb Z}}

\def\proof{\noindent{\bf Proof.}\ \ }
\def\qed{\hfill $\Box$\medskip}

\begin{document}

\title{Homogeneous Schr\"odinger operators on half-line}

\author{Laurent Bruneau} 
\address[Laurent Bruneau]{University of Cergy-Pontoise, Department
  of Mathematics and 
UMR 8088 CNRS, 95000 Cergy-Pontoise, France}
\email{laurent.bruneau@u-cergy.fr} 

\author{Jan Derezi\'nski}
\address[Jan Derezi\'nski]{Department of Mathematical Methods in
  Physics,    University of Warsaw, Ho\.{z}a 74, 00-682 Warszawa, Poland}
\email{Jan.Derezinski@fuw.edu.pl}

\author{Vladimir Georgescu} 
\address[Vladimir Georgescu]{CNRS and 
  University of Cergy-Pontoise 95000 Cergy-Pontoise, France}
\email{vlad@math.cnrs.fr} 

\date{\today}

\begin{abstract}

  The differential expression $L_m=-\partial_x^2+(m^2-1/4)x^{-2}$ 
  defines a self-adjoint operator $H_m$ on $L^2(0,\infty)$ in a
  natural way when $m^2\geq1$.  We study the dependence of $H_m$ on
  the parameter $m$, show that it has a unique holomorphic extension
  to the half-plane $\Re m >-1$, and analyze spectral and scattering
  properties of this family of operators.

\end{abstract}

\maketitle

\section{Introduction}\label{s:intro}

For $m\ge 1$ real the differential operator
$L_m=-\partial_x^2+(m^2-1/4)x^{-2}$ with domain
$C_\c^\infty=C_\c^\infty(0,\infty)$ is essentially self-adjoint and
we denote by $H_m$ its closure.  Let $U_\tau$ be the group of
dilations on $L^2$, that is $(U_\tau f)(x)=\e^{\tau/2} f(e^\tau x)$.
Then $H_m$ is clearly homogeneous of degree $-2$, i.e.  $U_\tau H_m
U^{-1}_\tau= \e^{-2\tau}H_m$. The following theorem summarizes the
main results of our paper.

\begin{theorem}\label{th:imain}
  There is a unique holomorphic family $\{H_m\}_{\Re m>-1}$ such
  that $H_m$ coincides with the previously defined operator if
  $m\ge1$. The operators $H_m$ are homogeneous of degree $-2$ and
  satisfy $H_m^*=H_{\bar{m}}$. In particular, $H_m$ is
  self-adjoint if $m$ is real.  The spectrum and the essential
  spectrum of $H_m$ are equal to $[0,\infty[$ for each $m$ with $\Re
  m>-1$. On the other hand, for non real $m$ the numerical range of
  $H_m$ depends on $m$ as follows:
\begin{compactenum}
 \item[i)] If $0<\arg m\leq\pi/2$, then
   $\Num(H_m)=\{z\ |\ 0\leq\arg  z\leq 2\arg m\}$,
 \item[ii)] If $-\pi/2\leq\arg m< 0$, then
   $\Num(H_m)=\{z\ |\ 2\arg m \leq \arg z \leq 0 \}$,
 \item[iii)] If $\pi/2 < |\arg m|<\pi$, then $\Num(H_m)=\C$.
\end{compactenum}
If $\Re m>-1$, $\Re k>-1$ and $\lambda\notin [0,\infty[$ then 
$(H_m-\lambda)^{-1}-(H_k-\lambda)^{-1}$ is a compact operator.
\label{main-result}
\end{theorem}

We note that if $0\le m<1$ the operator $L_m$ is not essentially
self-adjoint.  If $0< m<1$ this operator has exactly two distinct
homogeneous extensions which are precisely the operators $H_m$ and
$H_{-m}$ defined in the theorem: they are the Friedrichs and Krein
extension of $L_m$ respectively. Theorem \ref{th:imain} thus shows
that we can pass holomorphically from one extension to the
other. Note also that $L_0$ has exactly one homogeneous extension,
the operator $H_0$ which is at the same time the Friedrichs and the
Krein extension of $L_0$.
We obtain these results via a rather complete analysis of the
extensions (not necessarily self-adjoint) of the operator $L_m$ for
complex $m$.

We are not aware of a similar analysis of the  holomorphic family 
$\{H_m\}_{\Re m>-1}$ in the
literature. Most of the literature seems to restrict itself to the case of
real $m$ and self-adjoint $H_m$. A detailed study of the case $m>0$ can be
found in \cite{EK}. The fact that in this case the operator $H_m$ is
 the Friedrichs extension of $L_m$
 is of course well known. However, even the analysis
of the case
$-1\leq m\leq0$ seems to have been
 neglected in the literature.

We note that similar results concerning the holomorphic dependence
in the parameter $\alpha$ of the operator
$H_\alpha=(-\Delta+1)^{1/2}-\alpha/|x|$ have been obtained in \cite{Y}
by  different techniques.

Besides results described in Theorem \ref{main-result}, we prove a number of
other properties of the Hamiltonians $H_m$. Among other things, we
 treat the
scattering theory of the operators $H_m$ for real $m$, see Section
\ref{s:scm}.  We obtain explicit formulas for the corresponding
 wave  and scattering
operators. Essentially identical formulas in the closely related
context of
 the Aharonov-Bohm Hamiltonians 
 were obtained independently by Pankrashkin and
Richard in a recent paper \cite{PR} .

The   scattering theory for $H_m$ suggests a question, which we
were not able to answer. We pose this question as
 an interesting open problem in Remark
 \ref{open-problem}: can the
holomorphic family $\{\Re(m)>-1\}\ni 
m\mapsto H_m$ be extended to  the whole complex plane?
To understand why it is not easy to answer this question let us mention that
for $\Re(m)>-1$, the resolvent set is non-empty, being equal
to $\C\backslash[0,\infty[$. Therefore, to prove that $\{\Re(m)>-1\}
\ni m\mapsto H_m$ is a
    holomorphic family,
it suffices
    to show that its resolvent is holomorphic. However, one can show that if 
an extension of  this family to $\C$ exists, then
  for $\{m\ |\ \Re m=-1,-2,\dots,\ \ \Im
m\neq0\}$ the operator $H_m$ will have an empty resolvent set.
Therefore, on this set we cannot use the resolvent of $H_m$.

Let us describe the organization of the paper.
In Section \ref{s:pre} we recall
some facts concerning holomorphic families of closed operators and
make some general remarks on homogeneous operators and their
scattering theory in an abstract setting. Section \ref{s:foo} is
devoted to a detailed study of the first order homogeneous
differential operators. We obtain there several  results, 
which are then 
used in Section \ref{s:main}  containing our main
results.  In Section \ref{s:hank} we give
explicitly the spectral representation of $H_m$ for real $m$ and in
Section \ref{s:scm} we treat their scattering theory. In the first
 appendix we recall some technical results on second order
differential operators. Finally, as an application of Theorem
\ref{th:imain}, in the  second appendix we consider the Aharonov-Bohm
Hamiltonian $M_\lambda$ depending on the magnetic flux $\lambda$
and describe various holomorphic homogeneous rotationally symmetric
extensions of the family $\lambda\to
M_\lambda$. For a recent review on Aharonov-Bohm
Hamiltonians we refer to \cite{PR} and references therein.

To sum up, we believe that the operators $H_m$ are interesting for many
reasons.
\begin{itemize}
\item They have several  interesting physical applications, eg. they appear in
  the decomposition of the 
    Aharonov-Bohm   Hamiltonian.
\item They have rather subtle and rich properties, illustrating various
  concepts of the operator theory in Hilbert spaces (eg. the Friedrichs and
 Krein
  extensions, holomorphic families of closed operators). Surprisingly
 rich is also
  the theory of the first order homogeneous operators $A_\alpha$, that we
  develop in Sect. \ref{s:foo}, which is closely related to the theory of
  $H_m$. 
\item Essentially all basic objects related to $H_m$, such as their
  resolvents, spectral projections, wave and scattering operators, can be
  explicitly computed.
\item A number of
 nontrivial identities involving special functions find an
  appealing operator-theoretical interpretation in terms of the operators
  $H_m$. Eg. the Barnes identity (\ref{barnes}) leads to the formula for wave
  operators. Let us mention also the Weber-Schafheitlin
  identity \cite{KR}, which can be used to describe the distributional kernel
  of powers of $H_m$.
\end{itemize}

\medskip

{\bf Acknowledgement} J.D. would like to thank H.~Kalf for useful
discussions. His research was supported in part by the grant N N201 270135 of
the Polish Ministry of Science and Higher Education. The research of L.B. 
is supported by the ANR project HAM-MARK (ANR-09-BLAN-0098-01) of the French
Ministry of Research.

\section{Preliminaries } \label{s:pre}

\subsection{Notations} \label{ss:not}

For an operator $A$ we denote by $\cD(A)$ its domain, $\sp (A)$ its
spectrum, and $\rs(A)$ its resolvent set. We denote by $\Num A$ the
(closure of the) numerical range of an operator $A$, that is
$$
\Num A:=\overline{\{(f|Af)\ |\ f\in\cD(A), \ \|f\|=1\}}.
$$
If $H$ is a self-adjoint operator $H$ then $\cQ(H)$ will
denote its form domain: $\cQ(H)=\cD(|H|^{1/2})$.

We set $\R_+=\; ]0,\infty[$. We denote by $\one_+$ the
    characteristic function of the subset $\R_+$ of $\R$.

We write $L^2$ for the Hilbert space $L^2(\R_+)$. We abbreviate
$C_\c^\infty=C_\c^\infty(\R_+)$, $H^1=H^1(\R_+)$ and
$H_0^1=H_0^1(\R_+)$. Note that $H^1$ and $H_0^1$ are the form domains
of the Neumann and Dirichlet Laplacian respectively on $\R_+$. If
$-\infty\leq a<b\leq\infty$ we set $L^2(a,b)=L^2(]a,b[)$ and
    similarly for $C_c^\infty(a,b)$, etc.

Capital letters decorated with a tilde will denote operators acting
on distributions. For instance, let $\tilde Q$ and $\tilde P$ be the
position and momentum operators on $\R_+$, so that $(\tilde
Qf)(x)=xf(x)$ and $(\tilde Pf)(x)=-\i\partial_x f(x)$, acting in the
sense of distributions on $\R_+$. The operator $\tilde Q$ restricted
to an appropriate domain becomes a self-adjoint operator on $L^2$,
and then it will be denoted $Q$.  The operator $\tilde P$ has two
natural restrictions to closed operators on $L^2$, $P_\min$ with
domain $H_0^1$ and its extension $P_\max$ with domain $H^1$. We have
$(P_{\min})^*=P_{\max}$.

The differential operator $\tilde D:=\frac12(\tilde P\tilde Q+\tilde
Q\tilde P)=\tilde P\tilde Q+\i/2$, when considered as an operator in
$L^2$ with domain $C_\c^\infty$, is essentially self-adjoint and its
closure $D$ has domain equal to $\{f\in L^2 \mid PQf\in L^2\}$.  The
unitary group generated by $D$ is the group of dilations on $L^2$,
that is $(\e^{\i\tau D}f)(x)=\e^{\tau/2} f(e^\tau x)$.

We recall the simplest version of the Hardy estimate. 

\begin{proposition}\label{prop:hardy} For any  $f\in H_0^1$,
\[
\|P_\min f\|\geq \frac12\|Q^{-1}f\|.
\]
Hence, if $f\in H^1$ then $\tilde Q^{-1}f\in L^2$ if and only if $f\in H_0^1$.
\end{proposition} 

\proof For any $a\in\R$, as a quadratic form  on $\ccf$ we have
\begin{equation*}
0 \leq  ( \tilde P+\i a \tilde Q^{-1})^*( \tilde P+\i a \tilde Q^{-1})
 = \tilde  P^2+\i a[\tilde P, \tilde Q^{-1}]+a^2\tilde  Q^{-2}=
\tilde P^2+a(a-1) \tilde Q^{-2}.
\end{equation*}
Since $a(a-1)$ attains its minimum for $a=\frac12$, one gets
$\|\tilde Pf\|\geq \frac12\|\tilde Q^{-1}f\|$ for $f\in \ccf$. By the
dominated convergence theorem and Fatou lemma this inequality will
remain true for any $f\in H_0^1$.
\qed


\subsection{Holomorphic families of closed operators}\label{ss:hol}

In this subsection we recall the definition of a holomorphic family of
closed operators.

The definition (or actually a number of equivalent definitions) of a
{\em holomorphic family of bounded operators} is quite obvious and
does not need to be recalled. In the case of unbounded operators the
situation is more subtle.

Suppose that $\Theta$ is an open subset of $\C$, $\cH$ is a Banach
space, and $\Theta\ni z\mapsto H(z)$ is a function  whose values
are closed operators on $\cH$.  We say that this is a {\em
  holomorphic family of closed operators} if for each $z_0\in\Theta$
there exists a neighborhood $\Theta_0$ of $z_0$, a Banach space
$\cK$ and a holomorphic family of bounded operators $\Theta_0\ni
z\mapsto A(z)\in B(\cK,\cH)$ such that $\Ran A(z)=\cD(H(z))$ and
\begin{equation*}
\Theta_0\ni z\mapsto H(z)A(z)\in B(\cK,\cH)\label{holo2}
\end{equation*}
is a holomorphic family of bounded operators.

We have the following practical criterion:

\begin{theorem}\label{crit} Suppose that $\{H(z)\}_{z\in\Theta}$ is a
function whose values are closed operators on $\cH$. Suppose in
addition that for any $z\in\Theta$ the resolvent set of $H(z)$
is nonempty. Then $z\mapsto H(z)$ is a holomorphic family of
closed operators if and only if for any $z_0\in\Theta$ there exists
$\lambda\in\C$ and a neighborhood $\Theta_0$ of $z_0$ such that
$\lambda\in\rs( H(z))$ for $z\in\Theta_0$ and $z\mapsto
(H(z)-\lambda)^{-1}\in B(\cH)$ is holomorphic on $\Theta_0$.
\end{theorem}

The above theorem indicates that it is more difficult to study
holomorphic families of closed operators that for some values of the
complex parameter have an empty resolvent set.

To prove the analyticity of the resolvent, the following elementary
result is also useful
\begin{proposition}\label{prop:holom} Assume $\lambda\in \rs(H(z))$
  for $z\in\Theta_0$. If there exists a dense set of vectors $\cD$
  such that $z\mapsto\langle f,(H(z)-\lambda)^{-1}g\rangle$ is
  holomorphic on $\Theta_0$ for any $f,g\in\cD$ and if $z\mapsto
  (H(z)-\lambda)^{-1}\in B(\cH)$ is locally bounded on $\Theta_0$,
  then it is holomorphic on $\Theta_0$. 
\end{proposition}





\subsection{Homogeneous operators}  \label{ss:hops}

Some of the properties of homogeneous Schr\"odinger operators follow
by general arguments that do not depend on their precise structure. In
this and the next subsections we collect such arguments. These two
subsections can be skipped, since all the results that are given here
will be proven independently.

Let $U_\tau$ be a strongly continuous one-parameter group of unitary
operators on a Hilbert space $\cH$. Let $S$ be an operator on $\cH$
and $\nu$ a non zero real number. We say that $S$ is
\emph{homogeneous (of degree $\nu$)} if $U_\tau S
U_\tau^{-1}=\e^{\nu\tau}S$ for all real $\tau$. More explicitly this
means $U_\tau \cD(S)\subset \cD(S)$ and $U_\tau
SU_\tau^{-1}f=\e^{\nu\tau}Sf$ for all $f\in \cD(S)$ and all
$\tau$. In particular we get $U_\tau \cD(S)=\cD(S)$.

We are really interested only in the case $\cH=L^2$ and
$U_\tau=\e^{\i\tau D}$ the dilation group but it is convenient to
state some general facts in an abstract setting. Then, since we
assumed $\nu\neq0$, there is no loss of generality if we consider
only the case $\nu=1$ (the general case is reduced to this one by
working with the group $U_{\tau/\nu}$).

Let $T$ be a homogeneous operator. If $T$ is closable and densely
defined then $T^*$ is homogeneous too. If $S\subset T$ then $S$ is
homogeneous if and only if its domain is stable under the operators
$U_\tau$.

Let $S$ be a homogeneous closed hermitian (densely defined)
operator. We are interested in finding homogeneous self-adjoint
extensions $H$ of $S$. Since a self-adjoint extension satisfies
$S\subset H\subset S^*$ we see that we need to find subspaces
$\cE$ with $\cD(S)\subset \cE\subset \cD(S^*)$ such that $U_\tau
\cE\subset \cE$ for all $\tau$. Such subspaces will be called
\emph{homogeneous}.

Set $\bra S^*f,g\ket-\bra f,S^*g\ket=\i\{f,g\}$. Then $\{f,g\}$ is a
hermitian continuous sesquilinear form on $\cD(S^*)$ which is zero
on $\cD(S)$ and a closed subspace $\cD(S)\subset \cE\subset
\cD(S^*)$ is the domain of a closed hermitian extension of $S$ if
and only if $\{f,g\}=0$ for $f,g\in \cE$. Such subspaces will be
called \emph{hermitian}.  Note the following obvious fact: for $f\in
\cD(S^*)$ we have $\{f,g\}=0$ for any $g\in \cD(S^*)$ if and only if
$f\in \cD(S)$.


If $T$ is a homogeneous operator and $\lambda\in\C$ is an eigenvalue
of $T$ then $\e^\tau\lambda$ is also an eigenvalue of $T$ for any
real $\tau$. In particular, a homogeneous self-adjoint operator
cannot have non-zero eigenvalues and its spectrum is $\R$, or
$\R_+$, or $-\R_+$, or $\{0\}$ (note that, since $U_\tau$ is a
strongly continuous one-parameter group, the least closed subspace
which contains an eigenvector and is stable under all the $U_\tau$
and all functions of the operator is separable).

The following result, due to von Neumann, is easy to prove:

\begin{proposition}\label{pr:ev}
Let $S$ be a positive hermitian operator with deficiency indices
$(n,n)$ for some finite $n\geq1$. Then for each $\lambda<0$ there is
a unique self-adjoint extension $T_\lambda$ of $S$ such that
$\lambda$ is an eigenvalue of multiplicity $n$ of
$T_\lambda$. Moreover, the negative spectrum of $T_\lambda$ is equal
to $\{\lambda\}$. In particular, if $S$ is homogeneous then
$T_\lambda$ is not homogeneous, so $S$ has non-homogeneous
self-adjoint extensions.
\end{proposition}

\proof It suffices to take
$\cD(T_\lambda)=\cD(S)+\ker(S^*-\lambda)$.\qed

Recall that the Friedrichs and Krein extensions of a positive
hermitian operator $S$ are positive self-adjoint extensions $F$ and
$K$ of $S$ uniquely defined by the following property: any positive
self-adjoint extension $H$ of $S$ satisfies $K\leq H \leq F$ (in the
sense of quadratic forms). Then a self-adjoint operator $H$ is a
positive self-adjoint extension of $S$ if and only if $K\leq H \leq
F$.

\begin{proposition}\label{re:KF}
If $S$ is as in Proposition \ref{pr:ev} and if the Friedrichs and
Krein extensions of $S$ coincide then any other self-adjoint
extension of $S$ has a strictly negative eigenvalue. 
\end{proposition}

\proof Indeed, such an extension will not be positive and its strictly negative spectrum
consists of eigenvalues of finite multiplicity.\qed

It is clear that any homogeneous positive hermitian operator has homogeneous self-adjoint extensions.


\begin{proposition}\label{pr:KF}
If $S$ is a homogeneous positive hermitian operator then the Friedrichs and Krein extensions of $S$ are homogeneous.
\end{proposition}
\proof For any $T$ we set $T_\tau=\e^{-\tau}U_\tau T U_{\tau}^{-1}$. Thus homogeneity means $T_\tau=T$. Then from $S\subset T \subset S^*$ we get $S\subset T_\tau \subset S^*$.  Clearly $F_\tau$ is a self-adjoint operator and is a positive extension of $S$ hence $F_\tau\leq F$. Then we also have
$F_{-\tau}\leq F$ or $\e^{\tau}U_{-\tau} F U_{-\tau}^{-1}\leq F$ hence $F\leq F_\tau$, i.e. $F=F_\tau$. Similarly $K=K_\tau$.  \qed


\subsection{Scattering theory for homogeneous operators } \label{s:sch}

In this subsection we continue with the abstract framework of Subsection \ref{ss:hops}.

We shall consider couples of self-adjoint operators $(A,H)$ such that $H$ is homogeneous with respect to the unitary group $U_\tau=\e^{\i\tau A}$ generated by $A$, i.e.  $U_\tau H U_\tau^{-1}=\e^\tau H$ for all real $\tau$. We the say that \emph{$H$ is a homogeneous Hamiltonian} (with respect to $A$).  This
can be \emph{formally} written as $[\i A,H]=H$.  It is clear that $H$ is homogeneous if and only if
$U_\tau\varphi(H)U_\tau^{-1}=\varphi(\e^\tau H)$ holds for all real $\tau$ and all bounded Borel functions $\varphi:\sigma(H)\to\C$. Also, it suffices that this be satisfied for only one function $\varphi$ which generates the algebra of bounded Borel functions on the spectrum of $H$, for example for just one continuous injective function. If we set $V_\sigma=\e^{\i\sigma H}$ then another way of writing the homogeneity condition is $U_\tau V_\sigma=V_{\e^\tau\sigma}U_\tau$ for all real $\tau,\sigma$.

We shall call $(A,H)$ a \emph{homogeneous Hamiltonian couple}. We say that this couple is \emph{irreducible} if there are no nontrivial closed subspaces of $\cH$ invariant under $A$ and $H$, or
if the Von Neumann algebra generated by $A$ and $H$ is $B(\cH)$. A direct sum (in a natural sense) of homogeneous couples is clearly a homogeneous couple. Below $H>0$ means that $H$ is positive and
injective and similarly for $H<0$.

\begin{proposition}\label{pr:ah}
A homogeneous Hamiltonian couple $(A,H)$ is unitarily equivalent to a direct sum of copies of homogeneous couples of the form $(P,\e^Q)$ or $(P,-\e^Q)$ or $(A_0,0)$ with $A_0$ an
arbitrary self-adjoint operator. If $H>0$ then only couples of the first form appear in the direct sum. A homogeneous Hamiltonian couple is irreducible if and only if it is unitarily equivalent to
one of the couples $(P,\e^Q)$ or $(P,-\e^Q)$ on $L^2(\R)$, or to some $(A_0,0)$ with $A_0$ a real number considered as operator on the Hilbert space $\C$. A homogeneous couple is irreducible if and only if one of the operators $A$ or $H$ has simple spectrum (i.e. the Von Neumann algebra generated by it is maximal abelian) and in this case both operators have simple spectrum.
\end{proposition}

\proof
By taking above $\varphi$ equal to the characteristic function of
the set $\R_+$ then $-\R_+$ and $\{0\}$ we see that the
closed subspaces $\cH_+,\cH_-,\cH_0$ defined by $H>0,H<0,H=0$
respectively are stable under $U_\tau$.  So we have a
direct sum decomposition $\cH=\cH_+\oplus\cH_-\oplus\cH_0$ which is
left invariant by $A$ and $H$, hence $A=A_+\oplus A_-\oplus A_0$ and
similarly for $H$, the operator $H_+$ being homogeneous with respect
to $A_+$ and so on. Since $H_0=0$ the operator $A_0$ can be
arbitrary. The reduction to $\cH_-$ is similar to the reduction to
$\cH_+$, it suffices to replace $H_-$ by $-H_-$.

Thus in order to understand the structure of an arbitrary
homogeneous Hamiltonian $H$ it suffices to consider the case when
$H>0$. If we set $S=\ln H$ then by taking $\varphi=\ln$ above we get
$U_\tau SU_\tau^{-1}=\tau +S$ for all real $\tau$, hence the couple
$(A,S)$ satisfies the canonical commutation relations and so we may
us the Stone-Von Neumann theorem: $\cH$ is a direct sum of subspaces
invariant under $A$ and $S$ and the restriction of this couple to
each subspace is unitarily equivalent to the couple $(P,Q)$ acting
in $L^2(\R)$. Since $H=\e^S$ we see that the restriction of
$(A,H)$ is unitarily equivalent to the couple $(P,\e^Q)$ acting in
$L^2(\R)$. 
\qed

\begin{remark}\label{re:irred}{\rm Thus an irreducible homogeneous couple with $H>0$ is unitarily
equivalent to the couple $(P,\e^Q)$ on $\cH=L^2(\R)$. A change of
variables gives also the unitary equivalence with the couple $(D,Q)$
acting in $L^2(\R_+)$, where $D=(PQ+QP)/2$.}
\end{remark}

In the next proposition \emph{we fix a self-adjoint operator $A$
  with simple spectrum on a Hilbert space $\cH$ and assume that
  there is at least a homogeneous operator $H$ with $H>0$}. Then the
spectrum of $A$ is purely absolutely continuous and equal to the
whole real line by the preceding results. Moreover, the spectrum of
$H$ is simple, purely absolutely continuous and equal to
$\R_+$. Homogeneity refers always to $A$.

\begin{proposition}\label{pr:wv}
Assume that $H_1,H_2$ are homogeneous hamiltonians such that $H_k>0$. Then there is a Borel function $\theta:\R\to\C$ with $|\theta(x)|=1$ for all $x$ such that $H_2=\theta(A)H_1\theta(A)^{-1}$. If $\theta'$ is a second function with the same properties then there is $\lambda\in\C$ such that
$|\lambda|=1$ and $\theta'(x)=\lambda\theta(x)$ almost everywhere. If the wave operator $\Omega_+=\slim_{t\to+\infty}\e^{\i tH_2}\e^{-\i tH_1}$ exists then there is a function $\theta$ as
above such that $\Omega_+=\theta(A)$ and this function is uniquely determined almost everywhere. If the wave operator $\Omega_-=\slim_{t\to-\infty}\e^{\i tH_2}\e^{-\i tH_1}$ also exists then there is a uniquely determined complex number $\xi$ such that $\xi\Omega_-=\Omega_+$. In particular, the scattering
matrix given by $S=\Omega_-^*\Omega_+=\xi$ is independent of the energy. 
\end{proposition}

\proof As explained above the couples $(A,H_1)$ and $(A,H_2)$
are unitarily equivalent, hence there is a unitary operator $V$ on
$\cH$ such that $VAV^{-1}=A$ and $VH_1V^{-1}=H_2$. The
spectrum of $A$ is simple and $V$ commutes with $A$ so there is a
function $\theta$ as in the statement of the proposition such that
$V=\theta(A)$. If $W$ is another unitary operator with the same
properties as $V$ then $WV^{-1}$ commutes with $A$ and $H_2$. From
the irreducibility of $(A,H_2)$ it follows that $WV^{-1}$ is a
complex number of modulus one. Uniqueness almost everywhere is a
consequence of the fact that the spectrum of $A$ is purely
absolutely continuous and equal to $\R$.

Assume that $\Omega_+$ exists.  If we denote $\sigma=\e^{-\tau}$
then
$$
\e^{\i tH_2}\e^{-\i tH_1}U_\tau = \e^{\i tH_2}U_\tau\e^{-\i\sigma tH_1}=
U_\tau\e^{\i\sigma tH_2}\e^{-\i\sigma tH_1}
$$
hence $\Omega_+U_\tau=U_\tau\Omega_+$ for all real $\tau$. So the
isometric operator $\Omega_+$ belongs to the commutant $\{A\}'$, but
$\{A\}''$ is a maximal abelian algebra by hypothesis, so equal to
$\{A\}'$. Hence $\Omega_+$ must be a function $\theta(A)$ of $A$, in
particular it must be a normal operator, hence unitary. Now we
repeat the arguments above. Since the spectrum of $A$ is equal to
$\R$ and is purely absolutely continuous we see that $|\theta(x)|=1$
and is uniquely determined almost everywhere. Similarly, if
$\Omega_-$ exists then it is a unitary operator in $\{A\}''$. Thus
$S=\Omega_-^*\Omega_+$ is a unitary operator in $\{A\}''$ but also
has the property $H_1S=SH_1$. Since the couple $(A,H_1)$ is
irreducible we see that $S$ must be a number.  \qed


\section{Homogeneous first order  operators}\label{s:foo}

In this section we prove some technical results on homogeneous first
order differential operators which, besides their own interest, will
be needed later on.  

For each complex number $\alpha$ let $\widetilde A_\alpha$ be the
differential expression
\begin{equation}\label{eq:aalpha} 
\widetilde A_\alpha := \tilde P+\i\alpha \tilde Q^{-1}=
-\i\partial_x+\i\frac{\alpha}{x}= -\i x^{\alpha}\partial_x x^{-\alpha},
\end{equation}
acting on distributions on $\R_+$. Its restriction to $\ccf$ is a closable operator in $L^2$ whose closure will be denoted $\aam$. This is the \emph{minimal operator} associated to $\widetilde A_\alpha$. The \emph{maximal operator} $\aaM$ associated to $\widetilde A_\alpha$ is defined as the restriction of
$\widetilde A_\alpha$ to $\cD(A_\alpha^{\max}):=\{f\in L^2 \mid \widetilde A_\alpha f \in L^2 \}$.

The following properties of the operators $\aam$ and $\aaM$ are easy
to check:
\begin{compactenum}
\item[(i)] $\aam\subset\aaM$,
\item[(ii)] $(A_\alpha^{\min})^*=A_{-\overline{\alpha}}^{\max}$ and
$(A_\alpha^{\max})^*=A_{-\overline{\alpha}}^{\min}$,
\item[(iii)] $A_\alpha^{\min}$ and $A_\alpha^{\max}$ are homogeneous
  of degree $-1$.
\end{compactenum}

A more detailed description of the domains of the operators $\aam$ and $\aaM$ is the subject of the next proposition. We fix $\xi\in C^\infty_\c([0,+\infty)$ such that $\xi(x)=1$ for $x\leq 1$ and $\xi(x)=0$ for $x\geq 2$ and set $\xi_\alpha(x)=x^\alpha\xi(x)$.

\begin{proposition}\label{pr:mM}
\begin{compactenum}
\item[(i)] We have $\aam=\aaM$ if and only if $|\Re\alpha| \geq 1/2$. 
\item[(ii)] If $\Re\alpha\neq1/2$ then $\cD(\aam)=H_0^1$.
\item[(iii)] If $\Re\alpha=1/2$ then $H_0^1\subsetneq
  H^1_0+\C\xi_\alpha \subsetneq \cD(\aam)$ and $H_0^1$ is a core for
  $\aam=\aaM$. 
\item[(iv)] If $|\Re\alpha|<1/2$ then
  $\cD(\aaM)=H^1_0+\C\xi_\alpha$. In particular, if $|\Re\alpha|<1/2$
  and $|\Re\beta|< 1/2$ then $\cD(\aaM)\neq \cD(A_{\beta}^{\max})$. 
\end{compactenum}
\end{proposition}

Now we discuss the resolvent families. Let $\C_\pm=\{\lambda\in\C\mid \pm\Im\lambda>0\}$. The holomorphy of families of unbounded operators is discussed in \S\ref{ss:hol}. 

\begin{proposition}\label{pr:res}
(1) Let $\Re\alpha >-1/2$.  Then
\begin{compactenum}
\item[(i)]  $\rs(\aaM)=\C_-$.
\item[(ii)] If $\Im\lambda<0$ then the resolvent $(\aaM-\lambda)^{-1}$ is an integral operator with kernel
  \beq\label{eq:raaM}
    (\aaM-\lambda)^{-1}(x,y)=-\i \e^{\i \lambda(x-y)}\left(\frac{x}{y}\right)^\alpha \one_+(y-x).
  \eeq
\item[(iii)] The map $\alpha\mapsto\aaM$ is holomorphic in the region $\Re\alpha>-1/2$.
\item[(iv)] Each complex $\lambda$ with $\Im\lambda>0$ is a simple eigenvalue of $\aaM$ with $x^\alpha\e^{\i\lambda x}$ as associated eigenfunction.
\end{compactenum}

(2)  Let $\Re\alpha < 1/2$. Then 
\begin{compactenum}
\item[(i)] $\rs(\aam)=\C_+$. 
\item[(ii)]  If $\Im\lambda>0$ then the resolvent $(\aam-\lambda)^{-1}$ is an integral operator with kernel
  \beq\label{eq:raam}
    (\aam-\lambda)^{-1}(x,y)=\i \e^{\i \lambda(x-y)}\left(\frac{x}{y}\right)^\alpha \one_+(x-y).
  \eeq
\item[(iii)] The map $\alpha\mapsto\aam$ is holomorphic in the region $\Re\alpha<1/2$.
\item[(iv)] The operator $\aam$ has no eigenvalues.
\end{compactenum}
\end{proposition}

In some cases $\aam$ and $\aaM$ are generators of semigroups. We define the generator of a semigroup $\{T_t\}_{t\geq0}$ such that formally $T_t=\e^{\i tA}$.  Note that in \eqref{eq:sgrp1} the function
$f$ is extended to $\R$ by the rule $f(y)=0$ if $y\leq0$.

\begin{proposition}\label{pr:sgrp}
If $\Re\alpha\geq 0$ then $\i A_\alpha^{\max}$ is the generator of a
$C^0$-semigroup of contractions
\begin{equation}\label{eq:sgrp}
(\e^{\i t\aaM}f)(x)=x^\alpha(x+t)^{-\alpha}f(x+t), \quad t\geq0,
\end{equation}
whereas if $\Re\alpha\leq0$ the operator $-\i\aam$ is the generator of
a $C^0$-semigroup of contractions
\begin{equation}\label{eq:sgrp1}
(\e^{-\i tA_\alpha^{\min}}f)(x)=x^\alpha(x-t)^{-\alpha}f(x-t), \quad t\geq0.
\end{equation}
The operators $\i\aaM$ for $-1/2<\Re\alpha<0$ and $-\i\aam$ for
$0<\Re\alpha<1/2$ are not generators of $C^0$-semigroups of bounded
operators.
\end{proposition}

The remaining part of this section is devoted to the proof of the
these propositions. We begin with a preliminary fact.

\begin{lemma}\label{lm:rs} 
If $R$ and $S$ are closed operators such that $0\in \rs(R)$ then the operator $RS$ defined on the domain $\cD(RS):=\{ f\in \cD(S)\mid Sf\in \cD(R)\}$ is closed.
\end{lemma}

\proof Let $u_n\in\cD(RS)$ such that $u_n\to u$ and $RSu_n\to
v$. Then $u_n\in \cD(S)$ and $Su_n\in \cD(R)$, so that
$Su_n=R^{-1}RSu_n\to R^{-1}v$ because $R^{-1}$ is continuous. Since
$S$ is closed, we thus get that $u\in \cD(S)$ and
$Su=R^{-1}v$. Hence $Su\in \cD(R)$, i.e. $u\in \cD(RS)$, and
$v=RSu$.  \qed

Note that the Hardy estimate (Proposition \ref{prop:hardy}) gives
$\|\widetilde A_\alpha f\|\leq (1+2|\alpha|) \|Pf\|$ for all $f\in
H_0^1$. Since $\ccf$ is dense in $H_0^1$ we get $H_0^1\subset
\cD(\aam)$ for any $\alpha$. Our next purpose is to show that
$\cD(\aam)=H_0^1$ if $\Re\alpha\neq 1/2$, which is part (ii) of
Proposition \ref{pr:mM}. 

\begin{lemma}\label{lm:new}
  If $\Re\alpha\neq 1/2$ then $\cD(\aam)=H_0^1$.
\end{lemma}
\proof We set $\beta=\i(1/2-\alpha)$ and observe that it suffices to
prove that the restriction of $\widetilde A_\alpha$ to $H_0^1$ is a
closed operator in $L^2$ if $\Im\beta\neq0$. For this we shall use
Lemma \ref{lm:rs} with $R=D-\beta$ and $S$ equal to the self-adjoint
operator associated to $Q^{-1}$ in $L^2$. Then it suffices to show
that $\widetilde A_\alpha|_{H_0^1}=RS$.

The equality $\widetilde A_\alpha = (\widetilde D-\beta)Q^{-1}$,
where $\widetilde D=(PQ+QP)/2$ is the extension to distributions of
$D$, holds on the space of all distributions on $\R_+$, so we only
have to check that the domain of the product $RS$ is equal to
$H_0^1$ (because $\beta$ is not in the spectrum of the self-adjoint
operator $D$).  As discussed before, if $f\in H_0^1$ then
$Q^{-1}f\in L^2$, so $f\in\cD(S)$, and $PQQ^{-1}f=Pf\in L^2$, so
$Sf\in\cD(D)$. Thus $H_0^1\subset\cD(RS)$. Reciprocally, if
$f\in\cD(RS)$ then $f\in L^2$, $Q^{-1}f\in L^2$, and
$\widetilde{D}Q^{-1}f\in L^2$. But $\widetilde{D}Q^{-1}f\in L^2$ is
equivalent to $Pf\in L^2$, so $f\in H^1$. Since $Q^{-1}f\in L^2$ we
get $f\in H_0^1$.  
\qed

Our next step is the proof of part (1) of Proposition
\ref{pr:res}. Assume $\Re\alpha>-\frac12$. The last assertion of
part (1) of Proposition \ref{pr:res} is obvious so $\sp(\aaM)$
contains the closure of the upper half plane. We now show that if
$\Im\lambda<0$ then $\lambda\in\rs(\aaM)$ and the resolvent
$(\aaM-\lambda)^{-1}$ is an integral operator with kernel as in
\eqref{eq:raaM}.

The differential equation $(A_\alpha-\lambda)f=g$ is equivalent to
$\frac{\d}{\d x}(x^{-\alpha}\e^{-\i\lambda x}f(x))=\i
x^{-\alpha}\e^{-\i\lambda x}g(x)$. Assume $g\in L^2(0,\infty)$. We
look for a solution $f\in L^2(0,\infty)$ of the previous equation.
Since $\Im(\lambda)<0$ the function $x^{-\alpha}\e^{-\i\lambda
  x}g(x)$ is square integrable at infinity. We thus may define an
operator $R_\alpha^\max$ on $L^2$ by
$$
(R_\alpha^\max g)(x)=-\i \int_x^\infty \left(\frac{x}{y}\right)^\alpha
\e^{\i\lambda(x-y)}g(y)\d y, 
$$
i.e. $R_\alpha^\max$ is the integral operator with kernel given by \eqref{eq:raaM}.

\begin{lemma}\label{lem:amaxres1}
$R_\alpha^\max$ is a bounded operator in $L^2$.
\end{lemma}

\proof For shortness, we write $R$ for $R_\alpha^\max$. In the sequel we denote $\lambda=\mu+\i\nu$ and
$a=\Re\alpha$. By our assumptions, we have $\nu<0$ and $a>-1/2$. If $a\geq 0$ then the proof of the lemma is particularly easy because
\begin{equation*}
\int_0^\infty |R(x,y)|\d y = x^a\e^{-\nu x} \int_x^\infty
y^{-a}\e^{\nu y}\d y \leq \e^{-\nu x}\int_x^\infty \e^{\nu y}\d y = -\nu^{-1},
\end{equation*}
and similarly $\int_0^\infty |R(x,y)|\d x\leq-\nu^{-1}$. Then the
boundedness of $R$ follows from the Schur criterion. To treat the
case $-1/2<a<0$ we split the integral operator $R$ in two parts
$R_0$ and $R_1$ with kernels
\begin{equation*}
 R_0(x,y)=\one_{]0,1[}(x)R(x,y),\quad R_1(x,y)=\one_{[1,\infty[}(x)R(x,y). 
\end{equation*}
We shall prove that $R_1$ is bounded and $R_0$ is Hilbert-Schmidt.
For $R_1$ we use again the Schur criterion. If $x<1$ then
$\int_0^\infty |R_1(x,y)|\d y=0$ while if $x\geq1$ then
\[
\int_0^\infty |R_1(x,y)|\d y= x^a\e^{-\nu x} \int_x^\infty
y^{-a}\e^{\nu y}\d y.
\] 
We then integrate by parts twice to get
\begin{equation}\label{eq:part}
\int_0^\infty |R_1(x,y)|\d y=-\nu^{-1}-\frac{a}{\nu^2 x}
+\frac{a(a+1)}{\nu^2}x^a\e^{-\nu x}\int_x^\infty \e^{\nu y}y^{-a-2}\d y.
\end{equation}
Then, using $a>-1/2$, we estimate 
\[ 
x^a\e^{-\nu x}\int_x^\infty \e^{\nu y}y^{-a-2}\d y \leq x^a\int_x^\infty y^{-a-2}\d y=\frac{1}{(a+1)x},
\]
which, together with \eqref{eq:part}, proves that $\sup_{x\geq 1} \int_0^\infty |R_1(x,y)|\d y<+\infty$. Similarly $\int_0^\infty |R_1(x,y)|\d x=0$ if $y<1$ and for $y\geq1$
\[
\int_0^\infty |R_1(x,y)|\d x= y^{-a}\e^{\nu y} \int_1^y x^{a}\e^{-\nu x}\d y
\]
is estimated similarly.  We now prove that the operator $R_0$ is
Hilbert-Schmidt. We have
$$
\int_0^\infty \d x\int_0^\infty \d y |R_0(x,y)|^2 = \int_0^1\d x\,
x^{2a}\e^{-2\nu x} \int_x^\infty \d y\, y^{-2a}\e^{2\nu y}.
$$
Since $a$ and $\nu$ are strictly negative the integral $\int_0^\infty y^{-2a}\e^{2\nu y}\d y$ converges. Hence
$$
\int_0^\infty \d x\int_0^\infty \d y |R_0(x,y)|^2\leq C \int_0^1
x^{2a}\e^{-2\nu x}\d x, 
$$
which is convergent because $a>-1/2$.
\qed

So we proved that for $\Im(\lambda)<0$ the operator $R$ defines a bounded operator on $L^2$ such that $(\widetilde A_\alpha-\lambda)Rg=g$ for all $g\in L^2$. Hence, $R:L^2\to \cD(\aaM)$ and $(\aaM-\lambda)R=\one_{L^2}$.

Reciprocally, let $f\in \cD(\aaM)$ and set $g:=(\aaM-\lambda)f\in L^2$. The preceeding argument shows that $(\aaM-\lambda)(f-Rg)=0$. But $\aaM-\lambda$ is injective. Indeed,
if $(\aaM-\lambda)h=0$, then there exists $C\in\C$ such that
$h(x)=Cx^\alpha\e^{\i\lambda x}$ which is not in $L^2$ near infinity
unless $C=0$ (recall that $\Im(\lambda<0)$).

We have therefore proven that each $\lambda\in\C_-$ belongs to the
resolvent set of $\aaM$ and that $(\aaM-\lambda)^{-1}=R$. If we fix
such a $\lambda$ and look at $R=R(\alpha)$ as an operator valued
function of $\alpha$ defined for $\Re\alpha>-1/2$ then from the
preceding estimates on the kernel of $R$ it follows that
$\|R(\alpha)\|$ is a locally bounded function of $\alpha$. On the
other hand, it is clear that if $f,g\in\ccf$ then
$\alpha\mapsto\langle f,R(\alpha)g\rangle$ is a holomorphic
function. Thus, by Proposition \ref{prop:holom},
$\alpha\mapsto(\aaM-\lambda)^{-1}$ is holomorphic on
$\Re\alpha>-1/2$. This finishes the proof of point (1) of
Proposition \ref{pr:res}. The second part of the proposition
follows from the first part by using the relation
$\aam=(A_{-\overline\alpha}^{\max})^*$.

We now complete the proof of Proposition \ref{pr:mM} and consider
first the most difficult case when $\Re(\alpha)=1/2$.  The function
$\xi_\alpha$ is of class $C^\infty$ on $\R_+$, is equal to zero on
$x>2$, we have $\xi_\alpha\in L^2$, and $\widetilde A_\alpha
\xi_\alpha=0$ on $x<1$.  Hence $\xi_\alpha\in \cD(\aaM)$. On
the other hand $\xi_\alpha'\notin L^2$ (it is not square integrable
at the origin) so $\xi_\alpha\notin H_0^1$.

\begin{lemma} Let $\Re(\alpha)\geq1/2$. Then $\xi_\alpha\in \cD(\aam)$.
\end{lemma}

\proof The case $\Re\alpha>1/2$ is obvious since $\xi_\alpha\in
H^1_0$.  Now for $\Re\alpha=1/2$ we prove that $\xi_\alpha$ belongs
to the closure of $H_0^1$ in $\cD(\aaM)$ which is precisely
$\cD(\aam)$.  For $0<\varepsilon<1/2$ we define
$\xi_{\alpha,\varepsilon}$ as
\begin{equation*}
\xi_{\alpha,\varepsilon}(x) = \left\{\begin{array}{lcl}
    \frac{x}{\varepsilon}x^\alpha & 
\text{if} & x<\varepsilon, \\ \xi_\alpha(x) & \text{if} & x\geq
\varepsilon. \end{array} \right. 
\end{equation*}
For $x<\varepsilon$ one has
$\xi_{\alpha,\varepsilon}'(x)=\frac{\alpha+1}{\varepsilon}x^\alpha$.
Hence $\xi_{\alpha,\varepsilon}'\in L^2$ so that
$\xi_{\alpha,\varepsilon}\in H_0^1$. Moreover
$\|\xi_{\alpha,\varepsilon}-\xi_{\alpha}\|_{L^2}\to 0$ as
$\varepsilon\to 0$.  We then have
\begin{equation*}
\widetilde A_\alpha\xi_{\alpha,\varepsilon}(x) = \left\{\begin{array}{lcl}
-\frac{\i}{\varepsilon}x^\alpha & \text{if} & x<\varepsilon \\ 0 & \text{if} & \varepsilon\leq x<1 \end{array}\right. \quad \text{and} \quad  \widetilde A_\alpha\xi_{\alpha}(x)=0 \quad\text{if } x<1,
\end{equation*}
while $\widetilde A_\alpha\xi_{\alpha,\varepsilon}(x)=
\widetilde A_\alpha \xi_{\alpha}(x)$ if  $x\geq1$. Therefore
\begin{equation*}
 \|\widetilde A_\alpha\xi_{\alpha,\varepsilon}\|_{L^2}^2 = 
\int_0^\varepsilon \left| \frac{x^\alpha}{\varepsilon}\right|^2\d x+
\|\widetilde A_\alpha\xi_\alpha\|_{L^2}^2 = \frac12+\|\widetilde A_\alpha \xi_{\alpha}\|_{L^2}^2.
\end{equation*}
Thus $\xi_{\alpha,\varepsilon}\to\xi_{\alpha}$ in $L^2$, 
$\xi_{\alpha,\varepsilon}\in H_0^1\subset\cD(\aaM)$, and there is $C>0$ such that $\|\widetilde
A_\alpha\xi_{\alpha,\varepsilon}\|_{L^2}\leq C$ for any $\varepsilon$. Since $\aaM$ is closed this proves that $\xi_\alpha$ belongs to the closure of $H_0^1$ in $\cD(\aaM)$, i.e. $\xi_\alpha\in\cD(\aam)$.

\begin{lemma} \label{lm:old}
Let $\Re(\alpha)\geq1/2$. Then $\cD(\aam)=\cD(\aaM)$.
\end{lemma}

Fix $\lambda\in\C$ such that $\Im(\lambda)<0$, e.g. $\lambda=-\i$, and let
$R=(\aaM+\i)^{-1}$. $R$ is continuous from $L^2$ onto $\cD(\aaM)$, hence $R(\ccf)$ is dense in $\cD(\aaM)$. Let now $g\in \ccf$ and $0<c<d<\infty$ such that $\text{supp\,} g\subset[c,d]$. Then for any $x<c$,
\begin{eqnarray*}
f(x) = (Rg)(x) & = & -\i x^\alpha \e^x\int_c^d y^{-\alpha}\e^{-y}g(y)\d y\\ 
 & \sim & Cx^\alpha+Cx^\alpha(\e^x-1)\ \sim \ Cx^\alpha+D x^{\alpha+1}
\end{eqnarray*}
as $x\to 0$. Hence $f\in \C\xi_\alpha +H_0^1$. Therefore
$R(C^\infty_c)\subset \C\xi_\alpha +H_0^1\subset \cD(\aam)$. Since
$R(C^\infty_c)$ is dense in $\cD(\aaM)$, the same is true for $\cD(\aam)$. But $\aam$ is a closed operator and so $\cD(\aam)=\cD(\aaM)$.
\qed

\begin{lemma} If $\Re \alpha=1/2$, then $\C\xi_\alpha +H_0^1\neq \cD(\aaM)$.
\end{lemma}

\proof Let $R$ be as above and let 
$g(y)=y^{-\bar{\alpha}}|\ln(y)|^{-\gamma}\one_{]0,\frac12[}(y)$ 
where $\gamma>1/2$. Then $g\in L^2$ hence $Rg\in \cD(\aaM)$. On the other hand, for
$x\leq1/2$ we have
$$
Rg(x) = -\i x^\alpha \e^x\int_x^{\frac12} \frac{\e^{-y}}{y|\ln(y)|^\gamma}\d y\sim Cx^\alpha
|\ln(x)|^{1-\gamma}
$$
as $x\to 0$. In particular, if $\gamma<1$ then $Rg\notin \C\xi_\alpha +H_0^1$.
\qed

All the assertions related to the case $\Re\alpha=1/2$ of
Proposition \ref{pr:mM} have been proved. Since 
\begin{equation}\label{eq:adj}
\aam=\aaM \Longrightarrow A_{-\bar\alpha}^{\min} = A_{-\bar\alpha}^{\max}
\end{equation}
holds for any $\alpha$, we get $\aam=\aaM$ and so $\cD(\aaM)=H_0^1$
if $\Re\alpha=-1/2$.  We now turn to the case $|\Re(\alpha)|>1/2$
and show $\cD(\aaM)=\cD(\aam)=H_0^1$. Due to \eqref{eq:adj} it
suffices to consider the case $\Re\alpha>1/2$ which is
precisely the statements of Lemmas \ref{lm:new} and \ref{lm:old}. 
Now we prove (iv) of Proposition \ref{pr:mM}.

\begin{lemma}
If $|\Re\alpha|<1/2$, then  $\C\xi_\alpha +H_0^1= \cD(\aaM)$.
\end{lemma}

\proof Clearly, $\xi_\alpha\notin H_0^1$. We easily show that $\xi_\alpha\in\cD(\aaM)$.

Once again, let $R=(\aaM+\i)^{-1}$ and let $f\in\cD(\aaM)$. There exists $g\in L^2$ such that $f=Rg$, or
\[
f(x) = -\i x^\alpha \e^x \int_x^\infty \e^{-y}y^{-\alpha}g(y)\d y.
\] 
We show that $f\in \C \xi_\alpha+H_0^1$. Clearly, only the behaviour at the
origin matters. For $x<1$ decompose $f$ as
\[
f(x) =  -\i x^\alpha \e^x \int_0^\infty \e^{-y}y^{-\alpha}g(y)\d y
+\i x^\alpha \e^x \int_0^x \e^{-y}y^{-\alpha}g(y)\d y
=:  f_0(x)+f_1(x).
\]
Note that the first integral makes sense because $|\Re(\alpha)|<1/2$
so $\e^{-y}y^{-\alpha}$ is square integrable. Clearly
\[
f_0(x)=C x^\alpha \e^x=Cx^\alpha+ Cx^\alpha(\e^x-1)\in 
\C\xi_\alpha +H_0^1
\]
near the origin. We then prove that $f_1\in H_0^1$ near the
origin. By construction, $(A_\alpha+\i)f_1=g\in L^2$, so if we prove
that $Q^{-1}f_1$ is in $L^2$ near the origin, we will get $f_1\in H^1$
near the origin and hence $f_1\in H_0^1$ near the origin.

For any $0<x<1$ we can estimate (with $a=\Re\alpha$ as before)
\begin{equation}\label{eq:l2estim}
\frac{1}{x}|f_1(x)| = \frac{1}{x}\left|\int_0^x \e^{x-y} \left(\frac{x}{y}\right)^\alpha g(y)\d y \right|  \leq   C \int_1^{+\infty} t^{a-2}|g(\frac{x}{t})|\d t. 
\end{equation}
For any $t\geq 1$ let $\tau_t$ be the map in $L^2$ defined by
$(\tau_tg)(x)=g(x/t)$ and let $T=\int_1^\infty t^{a-2}\tau_t \d
t$. We have $\|\tau_t\|_{L^2\to L^2}=\sqrt{t}$ hence $T$ is a
bounded operator on $L^2$ with $\|T\|\leq \int_1^\infty t^{a-3/2}\d
t$ which converges since $a<1/2$. Together with (\ref{eq:l2estim}),
this proves that $\frac{1}{x}f_1(x)$ is square integrable on
$]0,1[$. This completes the proof of Proposition \ref{pr:mM}. 

It remains to prove Proposition \ref{pr:sgrp}. Since this is just a
computation, we shall only sketch the argument. Note that it
suffices to consider the case of $\aaM$ because then we get the
result concerning $\aam$ by taking adjoints. Let us denote
$A_0^{\max}=P_{\max}$, so $P_{\max}$ is the restriction to the
Sobolev space $H^1$ of the operator $P$. It is well-known and easy
to check that $P_{\max}$ is the generator of the contraction
semigroup $(\e^{\i tP_{\max}}f)(x)=f(x+t)$ for $t\geq0$ and $f\in
L^2$. Now if we write \eqref{eq:aalpha} as $\widetilde A_\alpha=
Q^\alpha P Q^{-\alpha}$ then \eqref{eq:sgrp} is formally obvious
because it is equivalent to 
\begin{equation*}\label{eq:sgrpq}
\e^{\i t\aaM}=Q^{\alpha}\e^{\i tP_{\max}}Q^{-\alpha}. 
\end{equation*}
For a rigorous justification, we note that the right hand side here
or in \eqref{eq:sgrp} clearly defines a $C_0$-semigroup of
contractions if (and only if) $\Re\alpha\geq0$ and then a
straightforward computation shows that its generator is $\aaM$. One
may note that $\ccf+\C\xi_\alpha$ is a core for $\aaM$ for all such
$\alpha$.


\section{Homogeneous second order operators}
\label{s:main}

\subsection{Formal operators}

For an arbitrary complex number $m$ we introduce the differential
expression
\begin{equation}\label{eq:hmdiff}
  \tilde L_m=\tilde P^2+(m^2-1/4)\tilde Q^{-2} =
  -\partial_x^2+\frac{m^2-1/4}{x^2} 
\end{equation}
acting on distributions on $\R_+$. Let $L_m^\m$ and $L_m^\M$ be the
minimal and maximal operators associated to it in $L^2$ (see
Appendix \ref{s:so}). It is clear that they are homogeneous
operators (of degree $-2$, we shall not specify this anymore). The
operator $L_m^\m$ is hermitian if and only if $m^2$ is a real
number, i.e. $m$ is either real or purely imaginary, and then
$(L_m^\m)^*=L_m^\M$. In general we have
\begin{equation*}\label{eq:ladjoint}
(L_m^\m)^*=L_{\bar m}^\M.
\end{equation*}
Note that \eqref{eq:hmdiff} does not make any difference between $m$
and $-m$. We will however see that $m$, not $m^2$, is the natural
parameter. In particular this will be clear in the construction of
other $L^2$ realizations of $L_m$, i.e. operators $H$ such that
$L_m^\m\subset H\subset L_m^\max$.

Observe also that one can factorize
$\tilde L_m$ as
\begin{equation}\label{eq:hmfactor}
\tilde L_m = \left(\tilde P+\i\frac{\bar m+\frac12}{\tilde Q}\right)^*
\left(\tilde P+\i\frac{m+\frac12}{\tilde Q}\right) = 
\widetilde A_{\bar m+\frac12}^*\widetilde A_{m+\frac12}
\end{equation}
where $\widetilde A_{\bar m+\frac12}^*$ is the formal adjoint of the differential expression $\widetilde A_{\bar m+\frac12}$.  The above expression makes a priori a difference between $m$ and $-m$, since $\tilde L_m$ does not depend on the sign of $m$ whereas the factorizations corresponding to $m$ and $-m$ are different. These factorizations provide one of the methods to distinguish between the various homogeneous extensions of $L_m^\m$. However, as we have seen in the previous section, one has to be carefull in the choice of the realization of $\widetilde A_{m+\frac12}$.


\subsection{Homogeneous holomorphic family}
\label{ssec:holom}

If $m$ is a complex number we set
\begin{equation}\label{eq:um}
\zeta_m(x)=x^{1/2+m}  \text{ if }  m\neq0 \quad {\rm and}  \quad
\zeta_{0}(x) \equiv \zeta_{+0}(x)=\sqrt{x}, \quad \zeta_{-0}(x)=\sqrt{x}\ln{x}.
\end{equation}
The notation is chosen in such a way that for any $m$ the functions $\zeta_{\pm m}$
are linearly independent solutions of the equation  $L_m u=0$.
Note that $\zeta_{\pm m}$ are both square integrable at the origin if and only if $|\Re m|<1$. 

We also choose $\xi\in C^\infty(\R_+)$ such that $\xi=1$ on $[0,1]$ and $0$ on $[2,\infty[$.

\begin{definition}
For $\Re (m)>-1$, we define $H_m$ to be the operator $L_m^\max$ restricted to $\cD(L_m^\min)+\C\xi \zeta_m$.
\end{definition}

Clearly, $H_m$ does not depend on the choice of $\xi$. Our first result concerning the family of operators $H_m$ is its analyticity with respect to the parameter $m$.

\begin{theorem}\label{th:holo}
$\{H_m\}_{\Re m>-1}$ is a holomorphic family of operators. More
  precisely, the number $-1$ belongs to the resolvent set of $H_m$ for
  any such $m$ and $m\mapsto(H_m+1)^{-1}\in\cB(L^2)$ is a holomorpic
  map.
\end{theorem}

Before we prove the above theorem, let us analyze the eigenvalue problem for $\tilde L_m$. The latter is closely related to Bessel's equation. In the sequel, $J_m$ will denote the Bessel functions of the first kind, i.e.
\begin{equation}\label{eq:bessel}
J_m(x):= \sum_{j=0}^\infty \frac{(-1)^j(x/2)^{2j+m}}{j! \Gamma(j+m+1)},
\end{equation}
and $I_m$ and $K_m$ the modified Bessel functions \cite{W}
\begin{equation}\label{eq:modifbessel}
I_m(x)=\i^{-m}J_m(\i x), \qquad K_m(x)=\frac{\pi}{2}\frac{I_{-m}(x)-I_m(x)}{\sin(m\pi)}.
\end{equation}

\bel\label{prop:resolvhm} For any $m$ such that $\Re(m)>-1$, the functions $\sqrt{x}I_m(x),\sqrt{x}K_m(x)$ form a basis of solutions of the differential equation
$-\partial^2_xu+(m^2-\frac14)\frac{1}{x^2}u=-u$ such that $\sqrt{x}I_m(x)\in L^2(]0,1[)$ and $\sqrt{x}K_m(x)\in L^2(]1,+\infty[)$. Besides, the Wronskian of these two solutions equals $1$.
\eel
\proof
If we introduce $w(x)=x^{-1/2}v(x)$, then $v$ satisfies $\widetilde L_m v=-v$ iff $w$ satisfies
$$
x^2w''(x)+xw(x)-(x^2+m^2)w=0,
$$
which is modified Bessel's differential equation. Linearly independent solutions of this equation are 
$(I_m,K_m)$. Therefore, a basis of solution for the equation $\widetilde
L_mu=-u$ is $(\sqrt{x}I_m(x),\sqrt{x}K_m(x))=:(u_0,u_\infty)$.

One has $I_m'(x)K_m(x)-I_m(x)K_m'(x)=-\frac{1}{x}$ (see \cite{W}),
and hence $W=u_0'u_\infty-u_0u_\infty'=1$. Moreover, $I_m(x)\sim
\frac{1}{\Gamma(m+1)}\left(\frac{x}{2}\right)^m$ as $x\to 0$
\cite{W}. Therefore, $u_0(x)$ is square integrable near the origin
iff $\Re(m)>-1$. On the other hand, $K_m(x)\sim
\sqrt{\frac{\pi}{2x}}\e^{-x}$ as $x\to\infty$, so that $u_\infty$ is
always square integrable near $\infty$. \qed


Note that $\sqrt{x} I_m(x)$ belongs to the domain of $H_m$ for all
$\Re(m)>-1$. Therefore, the candidate for the inverse of the operator
$H_m+1$ has kernel (cf. Proposition \ref{lm:green})
\[
G_m(x,y)=\left\{\begin{array}{lcc} \sqrt{xy}I_m(x)K_m(y) & {\rm if} &
x<y, \\ \sqrt{xy}I_m(y)K_m(x) & {\rm if} & x>y.  \end{array}\right.
\] 
We still need to prove that $G_m$ is bounded, which will be proven in the next lemma.

\begin{lemma} The map $m\mapsto G_m$ is a holomorphic family of bounded operators and it does not have 
a holomorphic extension to a larger subset of the complex plane.
\end{lemma}

\proof We prove that $G_m$ is locally bounded and that $m\mapsto\langle f,G_mg\rangle$ is analytic for $f,g$ in a dense set of $L^2$, so that the result follows from Proposition \ref{prop:holom}.

The modified Bessel functions depend analytically in $m$. Therefore the Green function $G_m(x,y)$ is an analytic function of the parameter $m$, and it is easy to see that for any $f,g\in C_\c^\infty(]0,+\infty[)$, the quantity $\langle f,(H_m+1)^{-1}g \rangle = \int \bar{f}(x)G_m(x,y)g(y) \d x\d y$ is analytic in $m$. Since $C_\c^\infty(]0,+\infty[)$ is dense in $L^2(0,+\infty)$ it remains to prove that $(H_m+1)^{-1}$ is locally bounded in $m$.

We shall split this resolvent as $G_m=G_m^{--}+G_m^{-+}+G_m^{+-}+G_m^{++}$ where $G_m^{\pm\pm}$ is the operator which has kernel $G_m^{\pm\pm}(x,y)$ with 
\begin{eqnarray*}
G_m^{--}(x,y) & = & G(x,y)\one_{]0,1]}(x)\one_{]0,1](y)},\\
G_m^{-+}(x,y) & = & G(x,y)\one_{]0,1]}(x)\one_{]1,\infty[(y)},\\
G_m^{+-}(x,y) & = & G(x,y)\one_{]1,\infty[}(x)\one_{]0,1](y)},\\
G_m^{++}(x,y) & = & G(x,y)\one_{]1,\infty[}(x)\one_{]1,\infty[(y)}.
\end{eqnarray*}
We control the norm of $G_m^{++}$ using Schur's Theorem (see \cite{Z}), whereas for the other terms, we estimate the $L^2$ norm of the kernel (this means in particular that $G_m^{--}$, $G_m^{-+}$ and $G_m^{+-}$ are actually Hilbert-Schmidt).

For that purpose, we use the explicit expression given in
Lemma \ref{prop:resolvhm} together with the following
estimates on the modified Bessel functions (see e.g. \cite{W})
\begin{itemize}
 \item as $x\to 0$
\begin{eqnarray}
I_m(x) & \sim & \frac{1}{\Gamma(m+1)}
\left(\frac{x}{2}\right)^m,\qquad m\neq -1,-2,\dots;\label{eq:besselfunc}\\
K_m(x) & \sim & \left\{ \begin{array}{lcl} \Re \left(  \Gamma(m) \left(\frac{2}{x}\right)^{m}\right)
  & {\rm if} &  \Re m=0,\ m\neq0;\\[1ex] -\ln\left(\frac{x}{2}\right)-\gamma 
  & {\rm if} & m=0,\\[1ex] \frac{\Gamma(m)}{2} \left(\frac{2}{x}\right)^{m} &   {\rm if} & \Re m>0;  \\ [1ex] \frac{\Gamma(-m)}{2} \left(\frac{x}{2}\right)^{m} & {\rm if} & \Re m<0. \end{array} \right.  \label{eq:besselfunc2}
\end{eqnarray}
 \item as $x\to \infty$
\begin{eqnarray}
I_m(x) & \sim & \frac{1}{\sqrt{2\pi x}}\e^x,\label{eq:besselfunc3}\\
K_m(x) & \sim & \sqrt{\frac{\pi}{2x}}\e^{-x}.\label{eq:besselfunc4}
\end{eqnarray}
\end{itemize}
The various constants which appear in (\ref{eq:besselfunc})-(\ref{eq:besselfunc4}) are locally bounded in $m$, except $\Gamma(m)$ as $m$ goes to zero, so that we may estimate the $G_m^{\pm\pm}(x,y)$ by
\begin{eqnarray}
|G_m^{--}(x,y)| & \leq & C_m |\Gamma(m)|\left(x^{1/2-|\nu|}y^{1/2+\nu}\one_{0<y<x<1}(x,y)\right.\label{eq:green1}\\
 & & \qquad\qquad\qquad\qquad \left.+ x^{1/2+\nu}y^{1/2-|\nu|}\one_{0<x<y<1}(x,y)\right),\nonumber\\
|G_m^{+-}(x,y)| & \leq & C_m\e^{-x}y^{\nu+1/2}\one_{]1,\infty[}(x)\one_{]0,1](y)},\nonumber\\
|G_m^{-+}(x,y)| & \leq & C_m x^{\nu+1/2}\e^{-y}\one_{]0,1]}(x)\one_{]1,\infty[(y)},\nonumber\\
|G_m^{++}(x,y)| & \leq & C_m \e^{-|x-y|}\one_{]1,\infty[}(x)\one_{]1,\infty[(y)},\nonumber
\end{eqnarray}
where $\nu=\Re(m)$ and $C_m$ are constants which depend on $m$ but
are locally bounded in $m$. The only problem is when $m=0$ where we
shall replace (\ref{eq:green1}) by
\begin{equation}\label{eq:green1bis}
|G_0^{--}(x,y)| \leq C\left(
  y^{1/2}|\ln(x)|\one_{0<y<x<1}(x,y)+x^{1/2}|\ln(y)|\one_{0<x<y<1}(x,y)\right).  
\end{equation}
Note also that the constant appearing in (\ref{eq:green1}) blows up as $m$ goes to zero due to the factor $|\Gamma(m)|$.

Straightforward computation lead to the following bounds
\begin{eqnarray*}
\|G_m^{--}\|^2_{L^2} & \leq &
\frac{C_m|\Gamma(m)|}{(\nu+1)(4+2\nu-2|\nu|)}, \qquad m\neq 0,\\ 
\|G_m^{-+}\|^2_{L^2} & \leq & \frac{C_m}{4(1+\nu)},\\
\|G_m^{+-}\|^2_{L^2} & \leq & \frac{C_m}{4(1+\nu)},\\
\|G_m^{++}\|_{L_x^\infty(L^1_y)} & \leq & 2C_m,\\
\|G_m^{++}\|_{L_y^\infty(L^1_x)} & \leq & 2C_m.
\end{eqnarray*}
This proves that $G_m^{--}$, $G_m^{-+}$ and $G_m^{+-}$ are
Hilbert-Schmidt operators whose norm is locally bounded in $m$
(except maybe for $G_m^{--}$ near 0) and using Schur's Theorem
$G^{++}_m$ is bounded with $\|G_m^{++}\|\leq 2C(m)$.

It remains to prove that $G_m^{--}$ is locally bounded around $0$. To this end we use $|K_m(z)|<C\frac{|x^m-x^{-m}|}{|m|}$ and estimate the Hilbert-Schmidt norm, where we set $\nu:=\Re m$:
\begin{eqnarray*}
\int_{0<x<y<1}|G_m^{--}(x,y)|^2\d x\d y & \leq & \frac{C}{|m|^2} \int_{0<x<y<1}xy|x^m|^2|y^m-y^{-m}|^2\\
& \leq & \frac{C'}{|m|^2}\left(\frac{1}{4\nu+2}+\frac14-\frac{2}{2\nu+4}\right) =\frac{C'}{4(\nu+1)(\nu+2)}
\end{eqnarray*}
As a conclusion, $G_m$ is locally bounded in $m$ for all $m$ such that $\Re(m)>-1$.
\qed

This proves that for $\Re(m)>-1$ the number $-1$ belongs to the
resolvent set of $H_m$, we have $G_m=(H_m+1)^{-1}$, and $H_m$ is a
holomorphic family of operators, cf. Proposition
\ref{prop:holom}. This proves Theorem \ref{th:holo}.

The next theorem gives more properties of the operators $H_m$. The
main technical point is that the differences of the resolvents
$R_{m'}(\lambda)-R_{m''}(\lambda)$ are compact operators, where we set
$R_m(\lambda)=(H_m-\lambda)^{-1}$ for $\lambda$ in the resolvent set
of $H_m$. For the proof we need the following facts.

\begin{lemma}\label{lm:holo}
Let $\Omega$ be an open connected complex set, $X$ a Banach space, $Y$
a closed linear subspace of $X$, and $F:\Omega\to X$ a holomorphic
map.  If $F(z)\in Y$ for $z\in \omega$ where $\omega\subset\Omega$ has
an accumulation point in $\Omega$, then $F(z)\in Y$ for $z\in\Omega$.
\end{lemma}

Indeed, all derivatives of $F$ at an accumulation point of
$\omega$ in $\Omega$ can be computed in terms of $F|_{\omega}$,
 hence belong to the closed subspace generated
 the
$F(z)$ with $z\in \omega$.

\begin{lemma}\label{lm:relcomp}
Let $S,T$ be two closed operators on a Banach space $\cH$ and let
$K(\lambda)=(S-\lambda)^{-1}-(T-\lambda)^{-1}$. If $K(\lambda)$ is
compact for some $\lambda\in\rs(S)\cap\rs(T)$ then $K(\lambda)$ is
compact for all $\lambda\in\rs(S)\cap\rs(T)$.
\end{lemma}
\proof We denote $S_\lambda=(S-\lambda)^{-1}$ and
$S_{\lambda\mu}=(S-\lambda)(S-\mu)^{-1}$ and use similar notations
when $S$ is replaced by $T$. Then $S_\lambda=S_\mu S_{\mu\lambda}$
hence $K(\lambda)=K(\mu) S_{\mu\lambda} +T_\mu (S_{\mu\lambda} -
T_{\mu\lambda})$.  If $K(\mu)$ is compact then the first term on the
right hand side is compact. For the second term we note that
\[
S_{\mu\lambda} - T_{\mu\lambda} = 
S_{\lambda\mu}^{-1} - T_{\lambda\mu}^{-1}=
\left(1+(\mu-\lambda) S_\mu\right)^{-1} -
\left(1+(\mu-\lambda) T_\mu\right)^{-1} =
 (\mu-\lambda)  S_{\mu\lambda}K(\mu) T_{\mu\lambda}
\]
and the last expression is a compact operator.
\qed

\begin{theorem}\label{th:reso}
For any $\Re(m)>-1$ we have $\sp (H_m)=\bar\R_+$ and if
$\lambda\in\C\backslash\bar\R_+$ then $R_m(\lambda)-R_{1/2}(\lambda)$
is a compact operator. If $R_m(\lambda;x,y)$ is the integral kernel of
the operator $R_m(\lambda)$, then for $\Re k>0$ we have:
\begin{equation}\label{eq:reso} 
R_m(-k^2; x,y) = 
\left\{\begin{array}{lcc} \sqrt{xy}I_m(kx)K_m(ky) & {\rm if} & x<y, \\
\sqrt{xy}I_m(ky)K_m(kx) & {\rm if} & x>y.  \end{array}\right.
\end{equation}
\end{theorem}
\proof We first show that $G_m-G_{1/2}$ is compact for all $m$. From
Lemma \ref{lm:holo} it follows that it suffices to prove this for
$0<m<1/2$. In this case $H_m$ is a positive operator and we have
$H_m=H_{1/2}+V$ in the form sense, where $V(x)=a x^{-2}$  with
$a=m^2-1/4$ hence $-1/4<a<1/4$. The Hardy estimate
(Proposition \ref{prop:hardy}) implies $\pm V\leq 4|a| H_{1/2}$
and $4|a|<1$ so if we set $S=(H_{1/2}+\lambda)^{-1/2}$ with
 $\lambda>0$ we get
\[
\pm S V S \leq
4|a| H_{1/2} (H_{1/2}+\lambda)^{-1} \leq 4|a|<1.
\]
Thus $\|S V S \|<1$.  From 
$H_m+\lambda = S^{-1}(1+ SVS ) S^{-1}$ we obtain
\[
(H_m+\lambda)^{-1} = S (1+ SVS )^{-1} S = (H_{1/2}+\lambda)^{-1/2}
+ \sum_{n>0} (-1)^n S(SVS)^{n-1} SVS^2
\] 
where the series is norm convergent. Hence 
$R_m(-\lambda) - R_{1/2}(-\lambda)$ is compact if $SVS^2$ is compact
(recall that we assume $0<m<1/2$).

We now prove that $SVS^2$ is a compact operator. Note that
$S^2=(H_{1/2}+\lambda)^{-1}$ and $H_{1/2}$ is the Dirichlet Laplacian,
so that $S^2 L^2 =H_0^1\cap H^2$ and $S L^2 =H_0^1$. Thus we have to
show that $V$ when viewed as operator $H_0^1\cap H^2 \to H^{-1}$ is
compact. Clearly this operator is continuous, in fact $V$ is
continuous as operator $H_0^1 \to H^{-1}$. Moreover, $H_0^2$ is the
subspace of $H_0^1\cap H^2$ defined by $f'(0)=0$ hence is a closed
subspace of codimension one of $H_0^1\cap H^2$. Thus it suffices to
prove that $V:H_0^2\to H^{-1}$ is compact. Let $\theta$ be a
$C^\infty$ function which is equal to one on for $x<1$ and equal to
zero if $x>2$. Clearly $(1-\theta)V$ is a compact operator $H_0^2\to
L^2$ and so it suffices to prove that $\theta V:H_0^2\to H^{-1}$ is
compact. Again it is clear that $\theta : L^2 \to H^{-1}$ is compact,
so it suffices to show that  $V:H_0^2\to L^2$ is continuous. If
$f\in C^\infty_0$ then
\[
V(x)f(x)=x^2V(x) \int_0^x \frac{x-y}{x^2}  f''(y) \d y =
x^2V(x) \int_0^1 (1-t)  f''(tx) \d t.
\]
So if $c=\sup_x |x^2 V(x)|$ then
\[
\|Vf\|\leq c \int_0^1 (1-t) \|f''(t\cdot)\| \d t=
c \int_0^1 (1-t) t^{-1/2} \d t \|f''\| =
\frac{4c}{3} \|f''\|
\]
hence $V:H_0^2\to L^2$ is continuous.

Thus we proved that $R_m(-1)-R_{1/2}(-1)$ is a compact operator if
$\Re(m)>-1$. From Lemma \ref{lm:relcomp} it follows that
$R_m(\lambda)-R_{1/2}(\lambda)$ is compact if $\lambda$ is in the
resolvent set of $H_m$ and of $H_{1/2}$. We have $\sp
(H_{1/2})=\bar\R_+$ and we now show that $\sp (H_m)=\bar\R_+$.
Clearly the operator $G_{1/2}$ is self-adjoint, its spectrum is the
interval $[0,1]$, and we have $G_m=G_{1/2}+K$ for some compact
operator $K$. Thus if $z\not\in[0,1]$ we have
\[
G_m-z=(G_{1/2}-z) \left[ 1+(G_{1/2}-z)^{-1}K \right] \equiv
(G_{1/2}-z) \left[ 1+K(z) \right]
\]
where $K(\cdot)$ is a holomorphic compact operator valued function on
$\C\setminus [0,1]$ such that $\|K(z)\|\to 0$ as $z\to\infty$. From
the analytic Fredholm alternative it follows that there is a discrete
subset $N$ of $\C\setminus [0,1]$ such that $1+K(z)$ is a bijective
map $L^2\to L^2$ if $z\not\in [0,1] \cup N$. Thus $G_m-z$ is a
bijective map in $L^2$ if $z\not\in N \cup [0,1]$. The function 
$z\mapsto \lambda =z^{-1}-1$ is a homeomorphism of $\C\setminus\{0\}$
onto $\C\setminus\{-1\}$ which sends $]0,1]$ onto $\bar\R_+$ hence the
image of $N$ through it  is a set $M$  whose accumulation points
belong to $\bar\R_+\cup\{-1\}$. If 
$\lambda\not\in\bar\R_+\cup\{-1\}\cup M$ then
\[
(\lambda+1)^{-1}-(H_m+1)^{-1}=
(\lambda+1)^{-1} (H_m-\lambda) (H_m+1)^{-1}
\]
and the left hand side is a bijection in $L^2$ hence $H_m-\lambda$ is
a bijective map $\cD(H_m)\to L^2$ so $\lambda$ belongs to the
resolvent set of $H_m$.  Thus the spectrum of $H_m$ is included in
$\bar\R_+\cup\{-1\}\cup M$. But $H_m$ is homogeneous so $\sp (H_m)$
must be a union of half-lines. Since it is not empty, it has to be
equal to $\bar\R_+$.

The explicit form of the kernel of $R_m(\lambda)$ given in
\eqref{eq:reso} may be proven by a minor variation of the arguments of
the proof of Theorem \ref{th:holo} based on more refined estimates for
the modified Bessel functions. Since we shall not need this formula,
we do not give the details.  \qed

\begin{remark}\label{re:gen}{\rm
We describe here in more abstract terms the main fact behind the
preceding proof. Let $H_0$ be a self-adjoint operator on a Hilbert
space $\cH$ with form domain $\cK=\cD(|H_0|^{1/2})$ and let $V$ be a
continuous symmetric sesquilinear form on $\cK$. If $V$ when viewed
as operator $\cK\to\cK^*$ is compact then it is easy to prove that
the form sum $H=H_0+V$ is well defined and that
$(H-z)^{-1}-(H_0-z)^{-1}$ is a compact operator on $\cH$ (in fact,
also as operator $\cK^*\to\cK$). This compactness condition on $V$
is never satisfied if $H_0$ and $V$ are homogeneous of the same
orders so this criterion is useless in our context. But our argument
requires only that $V$ be compact as operator $\cD(H_0)\to\cK^*$ and
this property holds in the case of interest here.  }\end{remark}


\subsection{Domain of the minimal and maximal operator}

In this subsection we analyze the operators $L_m^\min$ and $L_m^\max$.

\begin{proposition} If $|\Re m|<1$ then $L_m^\m\subsetneq L_m^\M$ and $\cD(L_m^\m)$ is a closed subspace of codimension two of $\cD(L_m^\M)$.
\end{proposition}

\proof In this case, we have two solutions of $L_mu=0$ that are in $L^2$ around $0$. Hence, the result follows from Proposition \ref{pr:sade}.
\qed

\begin{proposition}\label{pr:m=M}
If $|\Re m|\geq1$ then $L_m^\m=L_m^\M$. Hence, for $\Re(m)\geq1$,
 $H_m=L_m^\min=L_m^\max$.
\end{proposition}

\proof We use the notation of the proof of Lemma \ref{prop:resolvhm}. We know that the operator $G_m$ is continuous in $L^2$, that the functions $u_0$ and $u_\infty$ are uniquely defined modulo constant
factors, and there are no solutions in $L^2$ of the equation $(\tilde L_m+1)u=0$. Lemma \ref{lm:green} says that $(\tilde L_m +1)G_mg=g$ for all $g\in L^2$, hence $(L_m^\max+1)G_m=1$ on $L^2$. In particular
$G_m:L^2\to\cD(L_m^\M)$ is continuous.  More explicitly, we have
\begin{equation*}\label{eq:ex}
(G_mg)(x)=u_0(x)\int_x^\infty u_\infty(y)g(y)\d y + u_\infty(x) \int_0^x u_0(y)g(y)\d y.
\end{equation*}
Now we shall use the following easily proven fact.

\textit{Let $E$ be a normed space and let $\varphi,\psi$ be linear
  functionals on $E$ such that a linear combination $a\varphi+b\psi$
  is not continuous unless it is zero. Then $\ker\varphi\cap\ker\psi$
  is dense in $E$. }

We take $E=C_c^\infty$ equipped with the $L^2$ norm and $\varphi(g)=\int_0^\infty u_0(x)g(x)\d x$, $\psi(g)=\int_0^\infty u_\infty(x)g(x)\d x$. The linear combination $a\varphi+b\psi$ is given by a similar expression with $u=au_0+bu_\infty$ as integrating function.  Since $(\tilde L_m+1)u=0$ we have $u\in L^2$ only if $u=0$. Thus $E_0=\ker\varphi\cap\ker\psi$ is dense in $E$. It is clear that $G_m
E_0\subset C_c^\infty$. Hence by continuity we get $G_mL^2\subset\cD(L_m^\m)$ and thus $(L_m^\m+1)G_m=1$ on $L^2$. On the other hand it is easy to show that $G_m(\tilde L_m+1)f=f$ if $f\in C_c^\infty$,
hence $G_m(L_m^\m+1)=1$ on $\cD(L_m^\m)$. Thus $L_m^\m+1:\cD(L_m^\m)\to L^2$ is a bijective map.  Since $L_m^\M+1$ is an extension of $L_m^\m+1$ and is injective, we must have $L_m^\m=L_m^\M$.
\qed

If $m=1/2$, then clearly  $\cD(L_m^\m)=H_0^2$. If $m\neq1/2$ then $\cD(L_m^\m)\neq H_0^2$. However, the functions from $\cD(L_m^\m)$ behave at zero as if they were in $H_0^2$ with the exception of the case $m=0$.

\begin{proposition}\label{pr:lm0} Let $f\in\cD(L_m^\m)$.
\begin{compactenum}
\item[(i)] If  $m\neq0$, then $f(x)=o(x^{3/2})$ and $f'(x)=o(x^{1/2})$ as $x\to0$.
\item[(ii)] If $m=0$, then $f(x)=o(x^{3/2} \ln x)$ and $f'(x)=o(x^{1/2} \ln x)$  as $x\to0$.
\item[(iii)] For any $m$, $\cD(L_m^\m)\subset H_0^1$. 
\end{compactenum}
\end{proposition}

\proof Since $\tilde L_m$ does not make any difference between $m$ and $-m$, we may assume $\Re m\geq0$.

Assume first $\Re m\geq 1$. If $f\in\cD(L_m^\m)$ and $g=(L_m^\m+1)f$, then $f=G_mg$ and hence
$f=u_0 g_\infty+u_\infty g_0$ and $f'=u_0'g_\infty-u_\infty'g_0$ with $g_0(x)=\int_0^x
u_0(y)g(y)\d y$ and $g_\infty(x)=\int_x^\infty u_\infty(y)g(y)\d y$.
The functions $u_0$ and $u_\infty$ are of Bessel type and their behaviour at zero is known, see
(\ref{eq:besselfunc2}). More precisely if we set $\mu=\Re m$ then we have 
\[
u_0(x)=O(x^{\mu+1/2}), \quad u_0'(x)=O(x^{\mu-1/2}), \quad
u_\infty(x)=O(x^{-\mu+1/2}), \quad u_\infty'(x)=O(x^{-\mu-1/2}).
\]
Then  for $x<1$ we have
\begin{eqnarray*}
|u_0(x) g_\infty(x)| & \leq & C x^{\mu+1/2} \left(\int_x^1 y^{-\mu+1/2}|g(y)| \d y +
  \int_1^\infty |u_\infty(y)g(y)| \d y \right) \\
 & \leq & C x^{\mu+1/2} \left( \left(\frac{x^{2-2\mu}-1}{2\mu-2}\right)^{1/2}
  +\|u_\infty\|_{L^2(1,\infty)} \right) \|g\|,
\end{eqnarray*}
which is $O(x^{3/2})$. We have $u_\infty g_0=o(x^{3/2})$ by a simpler argument. Let $F$ be the Banach space consisting of continuous functions on $I=\,]0,1[$ such that $\|h\|_F\equiv \sup_{x\in I}x^{-3/2}|h(x)|<\infty$ and for $g\in L^2$ let $Tg$ be the restriction of $G_mg$ to $I$. By what we
have shown we have $TL^2\subset F$,  hence, by the closed graph theorem, $T:L^2\to F$ is a continuous operator. With the notations of the proof of Proposition \ref{pr:m=M}, if $g\in E_0$ then $Tg$ is
equal to zero near zero, so $T$ sends the dense subspace $E_0$ of $L^2$ into the closed subspace $F_0$ of $F$ consisting of functions such that $x^{-3/2}h(x)\to0$ as $x\to0$. By continuity we get
$TL^2\subset F_0$, hence $f(x)=o(x^{3/2})$. A similar argument based on the representation $f'=u_0'g_\infty-u_\infty'g_0$ gives $f'(x)=o(x^{1/2})$.

We treat now the case $0\leq\Re m<1$. Now all the solutions of the
equation $L_mu=0$ are square integrable at the origin, hence we may 
use Proposition \ref{pr:domto} with $v_\pm$ proportional to 
$\zeta_{\pm m}$. A straightforward computation gives for $m\neq0$ 
\begin{equation*}
|v_+(x)|\|v_-\|_x+|v_-(x)|\|v_+\|_x \leq C x^{3/2},
\quad
|v_+'(x)|\|v_-\|_x+|v_-'(x)|\|v_+\|_x \leq C x^{1/2}
\end{equation*} 
while if $m=0$ then
\begin{equation*}
|v_+(x)|\|v_-\|_x+|v_-(x)|\|v_+\|_x \leq C x^{3/2} (|\ln x|+1),
\quad
|v_+'(x)|\|v_-\|_x+|v_-'(x)|\|v_+\|_x \leq C x^{1/2} (|\ln x|+1).
\end{equation*} 
This finishes the proof.
\qed

We describe now some consequences of the representations
\eqref{eq:sol} and \eqref{eq:sold} in the present context.  We say
that a function $h$ \emph{is in $\cD(L_m^\m)$ near the origin} if
for some (hence any) function $\xi\in C_c^\infty(\R)$ which is one
on a neighbourhood of the origin we have $\xi
h\in\cD(L_m^\m)$. Assume $|\Re m|<1$ and let $f\in\cD(L_m^\M)$. Then
there are constants $a,b$ and a function $f_0$ which is in
$\cD(L_m^\m)$ near the origin such that
\begin{eqnarray}
f(x) & = & ax^{1/2-m}+bx^{1/2+m} + f_0(x) \quad\text{if }
m\neq0, \label{eq:soL} \\ 
f(x) & = & ax^{1/2}\ln x+bx^{1/2} + f_0(x) \quad\text{if }
m=0. \label{eq:soL0} 
\end{eqnarray}
These relations give by differentiation representations of $f'$. It
is clear by Proposition \ref{pr:lm0} that $f_0$ decays more rapidly
at zero than the other two terms, in particular the constants $a,b$
and the function $f_0$ are uniquely determined by $f$. This allows
one to state assertions converse to that of Proposition
\ref{pr:lm0}, for example:


\bep\label{re:domt} We have the following characterization of the
domain of the minimal operator: 
\begin{eqnarray*}
0< \Re m \leq \mu \leq 1 & \Rightarrow \ \cD(L_m^\min) & = \ \{f\in\cD(L_m^\M) \mid f(x)=o(x^{\mu+1/2})\}\\
 & & = \{f\in\cD(L_m^\M) \mid f'(x)=o(x^{\mu-1/2})\},\\
0 \leq \Re m< \mu \leq 1 & \Rightarrow \ \cD(L_m^\min) & = \ \{f\in\cD(L_m^\M) \mid    f(x)=O(x^{\mu+1/2})\}\\
 & & =\{f\in\cD(L_m^\M) \mid f'(x)=O(x^{\mu-1/2})\}.
\end{eqnarray*}
\eep


\subsection{Strict extensions of $L_m^\min$}
\label{ssec:lminextension}

Now we study the closed extensions of   $L_m^\min$ for $|\Re m|<1$. The first result is a particular case of Proposition \ref{pr:sade}. We recall that by a \emph{strict extension} of $L_m^\m$
we mean an operator $H$ such that $L_m^\m\subsetneq H \subsetneq L_m^\M$.
We denote by $W_x(f,g):=f(x)g'(x)-f'(x)g(x)$ the Wronskian of two functions $f$ and $g$ at point $x$, and take $\xi$ as in Section \ref{s:foo}.

\begin{proposition}\label{pr:extl}
Assume that $|\Re m|<1$. Let $u$ be a non-zero solution of $\tilde L_mu=0$. Then $W_0(u,f)=\lim_{x\to0} W_x(u,f)$ exists for each $f\in\cD(L_m^\M)$ and the operator $L_m^u$ defined as the restriction of $L_m^\M$ to the set of $f\in\cD(L_m^\M)$ such that $W_0(u,f)=0$ is a strict extension of $L_m^\m$. Reciprocally, each strict extension of $L_m^\m$ is of the form $L_m^u$ for some non-zero solution $u$ of $\tilde L_mu=0$, which is uniquely defined modulo a constant factor. We have $\cD(L_m^u)=\cD(L_m^\m)+\C\xi u$.
\end{proposition}

We shall describe now the homogeneous strict extensions of $L_m^\m$.
The case $|\Re m|\geq 1$ is trivial because $L_m^\m=L_m^\M$ is homogeneous.

\bep\label{df:hm} If $-1<\Re m<1$, then $H_m$ is the restriction of
$L_m^\M$ to the subspace defined by
\begin{equation}\label{eq:hm}
\lim_{x\to0} x^{m+1/2} \left( f'(x)-\frac{m+1/2}{x}f(x) \right)=0.
\end{equation}
\eep

\proof Observe that
\[
W_x(\zeta_m,f)=x^{m+1/2}f'(x)-(m+1/2){x}^{m-1/2}f(x) = x^{m+1/2}\left(f'(x)-\frac{m+1/2}{x}f(x)\right),
\]
so the limit from the left hand side of \eqref{eq:hm} exits for all $f\in\cD(L_\M)$ if $|\Re m|<1$. Hence, with the notation of Proposition \ref{pr:extl} we have $H_m=L_m^{\zeta_m}$ where $\zeta_m$ is defined in \eqref{eq:um}.
\qed


\begin{proposition}\label{pr:homc}
If $|\Re m|< 1$ and $m\neq0$ then $L_m^\m$ has exactly two homogeneous strict extensions, namely the operators $H_{\pm m}$. If $m=0$ then the operator $H_0$ is the unique homogeneous strict
extension of $L_m^\m$.
\end{proposition}

\proof Thanks to Proposition \ref{pr:extl} it suffices to see when the extension $L_m^u$ is homogeneous. If $(U_tf)(x)=\e^{t/2}f(\e^{t}x)$ then it is clear that $L_m^u$ is homogeneous if and only if its domain is stable under the action of $U_t$ for each real $t$. We have
\begin{eqnarray*}
W_0(u,U_tf) & = & \lim_{x\to0} \left(u(x)\e^{t/2}\frac{d}{dx}f(\e^t x)-u'(x)\e^{t/2}f(\e^t x)\right)\\
 & = & \e^{t/2}\lim_{x\to0} \left(\e^{t}{u}(x)f'(\e^{t}x)- {u}'(x)f(\e^{t}x)\right) \\
 & = & \e^{3t/2}\lim_{x\to0} \left({u}(\e^{-t}x)f'(x)-\e^{-t}{u}'(\e^{-t}x)f(x)\right).
\end{eqnarray*}
Thus we obtain
\begin{equation*}\label{eq:ut}
W_0(u,U_tf)=\e^{2t}W_0(U_{-t}u,f).
\end{equation*}
Let $u_t=\e^{2t}U_{-t}u$.  From Proposition \ref{pr:extl} we see that $\cD(L_u)=\cD(L_{u_t})$ for all real $t$ if and only if $u_t$ is proportional to $u$ for all $t$, which means that the function $u$
is homogeneous. Thus it remains to see which are the homogeneous solutions of the equation $L_mu=0$. Clearly $u_{\pm m}$ are both homogeneous and only they are so if $m\neq0$ and if $m=0$ then only
$u_{+0}$ is homogeneous.
\qed

\bep\label{re:domh} For $\Re m>0$, we have the following alternative characterizations of the domain of $H_m$:
\begin{eqnarray*}
0 <\mu\leq \Re m <1 & \Rightarrow & \cD(H_m)=\{f\in\cD(L_m^\M) \mid f(x)=o(x^{-\mu+1/2})\},\\
0 \leq\mu< \Re m <1 & \Rightarrow & \cD(H_m)=\{f\in\cD(L_m^\M) \mid f(x)=O(x^{-\mu+1/2})\}.
\end{eqnarray*}
\eep

\proof We use Propositions \ref{pr:lm0},
and the representations \eqref{eq:soL} and \eqref{eq:soL0}. \qed


\subsection{The hermitian case}

We shall consider now the particular case when $L_m^\m$ is hermitian, i.e. $m^2$ is a real number. Everything follows immediately from the preceding propositions and from the last assertion of Proposition \ref{pr:sade}. If $m$ is real or $m=\i\mu$ with $\mu$ real it suffices to consider the cases $m\geq0$ and $\mu>0$ because $L_m^\m=L_{-m}^\m$.

\begin{proposition}\label{pr:1}
The operator $H_m=L_m^\m$ is self-adjoint and homogeneous for $m^2\geq1$. When $m^2<1$ the operator $L_m^\m$ has deficiency indices $(1,1)$ and therefore admits a one-parameter family of self-adjoint extensions.
\begin{enumerate}
\item If $0<m<1$ and $0\leq\theta<\pi$ let $u_\theta$ be the function on $\R_+$ defined by
  \begin{equation}\label{eq:utheta}
  u_\theta(x)=x^{1/2-m}\cos\theta +x^{1/2+m}\sin\theta.
  \end{equation}
  Then each self-adjoint extension of $L_m^\m$ is of the form $H_m^\theta=L_m^{u_\theta}$ for a unique $\theta$. There are exactly two homogeneous strict extensions, namely the self-adjoint operators $H_m=H_m^{\pi/2}$ and $H_{-m}=H_m^{0}$.
\item If $m=0$ and $0\leq\theta<\pi$ let $u_\theta$ be the function on $\R_+$ defined by
  \begin{equation}\label{eq:u0theta}
  u_\theta(x)=x^{1/2} \ln x\cos\theta +x^{1/2}\sin\theta.
  \end{equation}
  Then each self-adjoint extension of $L_0^\m$ is of the form $H_0^\theta=L_0^{u_\theta}$ for a unique $\theta$. The operator $L_0^\m$ has exactly one homogeneous strict extension: this is the self-adjoint operator $H_0=H_0^{\pi/2}$.
\item Let $m^2<0$ so that $m=\i\mu$ with $\mu>0$. For
  $0\leq\theta<\pi$ let $u_\theta$ be the function given by
  \begin{equation}\label{eq:uitheta}
    u_\theta(x)=x^{1/2} \cos(\mu\ln x)\cos\theta +x^{1/2} \sin(\mu\ln x)\sin\theta.
  \end{equation}
  Then each self-adjoint extension of $L_m^\m$ is of the form
  $H_m^\theta=L_m^{u_\theta}$ for a unique $\theta$. The operator
  $L_m^\m$ does not have homogeneous self-adjoint extensions but has
  two homogeneous strict extensions, namely the operators $H_m$ and
  $H_{-m}$.
\end{enumerate}
\end{proposition}

We shall now study the quadratic forms associated to the self-adjoint operators $H_m^\theta$ for $0<m<1$.

We recall that $A_{1/2+m}^\m=A_{1/2+m}^\M$ if $\Re m \geq0$ and $A_{1/2-m}^\m=A_{1/2-m}^\M$ if $\Re m\geq 1$, see Proposition \ref{pr:mM}. Let us abbreviate $A_\alpha=\aam=\aaM$ when the minimal
and maximal realizations of $\widetilde A_\alpha$ coincide.

Recall also that for $0<m<1$,
\begin{equation*}\label{eq:pm}
\cD(A_{1/2-m}^\max) = H_0^1+\C\xi \zeta_{-m}.
\end{equation*}
By Proposition \ref{pr:mM} the operator $A_{1/2-m}^\max$ is closed in $L^2$ and $H_0^1$ is a closed subspace of its domain (for the graph topology) because $A_{1/2-m}^\max\upharpoonright_{H_0^1}=A_{1/2-m}^\m$ is also a closed operator. Note
that for $f\in H_0^1$ we have $f(x)=o(\sqrt{x})$ because
\begin{equation*}\label{eq:h01}
|f(x)|\leq \int_0^x|f'(x)|\d x\leq\sqrt{x}\|f'\|_{L^2(0,x)}.
\end{equation*}
Thus $\xi \zeta_{-m}\notin H_0^1$ and the sum $H_0^1+\C\xi\zeta_m$ is a topological direct sum in $\cD(A_{1/2-m}^\max)$. Hence each $f\in\cD(A_{1/2-m}^\max)$ can be uniquely written as a sum $f=f_0+c\xi \zeta_{-m}$ and the map $f\mapsto c$ is a continuous linear form on  $\cD(A_{1/2-m}^\max)$. We shall denote $\varkappa_m$ this form and observe that
\begin{equation*}\label{eq:vkappa}
\varkappa_m(f)=\lim_{x\to0} x^{m-1/2}f(x), \quad  f\in\cD(A_{1/2-m}^\max).
\end{equation*}
Note also that from Proposition \ref{pr:mM} we get $(A_{1/2-m}^\max)^*=A_{m-1/2}^\m$  in particular
$\cD\left((A_{1/2-m}^\max)^*\right) =H_0^1$.


\begin{proposition}\label{pr:quadr}
Let $0<m<1$ and $0\leq\theta<\pi$.
\begin{enumerate}
\item If $\theta=\pi/2$, then $\cD(H_m^{\pi/2})$ is a dense subspace of $H_0^1$ and if                     $f\in\cD(H_m^{\pi/2})$ then
  \begin{equation}\label{eq:quadrp}
  \bra f , H_m^{\pi/2} f \ket =\|A_{1/2+m} f\|^2=\|A_{1/2-m}^\max f\|^2.
  \end{equation}
  Thus $\cQ(H_m^{\pi/2})=H_0^1$. Moreover, we have $H_m^{\pi/2}=(A_{1/2+m})^*A_{1/2+m} = (A_{1/2-m}^\min)^*A_{1/2-m}^\min$.
\item Assume $\theta\neq\pi/2$.  Then $\cD(H_m^\theta)$ is a dense subspace of $\cD(A_{1/2-m}^\max)$ and   for each $f\in\cD(H_m^\theta)$ we have 
  \begin{equation}\label{eq:quadrt}
  \bra f , H_m^\theta f \ket=\|A_{1/2-m}^\max f\|^2 +m\sin(2\theta)|\varkappa_m (f)|^2.
  \end{equation}
  Thus $\cQ(H_m^\theta)=\cD(A_{1/2-m}^\max)$ and the right hand side of \eqref{eq:quadrt} is equal to the quadratic form of $H_m^\theta$. 
\end{enumerate}
\end{proposition}

\proof From Proposition \ref{pr:lm0}, the definition of $H_m$ and \eqref{eq:utheta} we get 
\begin{equation*}\label{eq:hth0}
\cD(H_m^\theta)=\cD(L_m^\m) + \C \xi u_\theta \subset H_0^1 + \C \xi u_\theta
= H_0^1 + \C \cos\theta \,\xi  \zeta_{-m},
\end{equation*}
because $\xi \zeta_m\in H_0^1$ if $m>0$. But $C_c^\infty\subset\cD(L_m^\m)$ so $\cD(H_m^\theta)$ is a dense subspace of $H_0^1 + \C \cos\theta \, \xi \zeta_{-m}$. Thus if $\theta=\pi/2$ we get $\cD(H_m^{\pi/2})\subset H_0^1$ and if $\theta\neq\pi/2$ then $\cD(H_m^{\theta})\subset\cD(A_{1/2-m}^\max)$ densely in both cases.

The relation $\|A_{1/2-m}^\max f\|^2=\|A_{1/2+m} f\|^2$ for $f\in H_0^1$ holds because both terms are continuous on $H_0^1$ by Hardy inequality and they are equal to $\bra f , \tilde L_m f\ket$ if $f\in
C_c^\infty$.

It remains to establish \eqref{eq:quadrt}. Let $f=f_0+c\xi u_\theta$ with $f_0\in\cD(L_m^\m)$ and $c\in\C$. Then $A_{1/2-m}^\max f\in L^2$ and  
\[
H_m^\theta f=\tilde L_mf=\widetilde A_{1/2-m}^*A_{1/2-m}^\max f\in L^2
\]
due to \eqref{eq:hmfactor}. Denote $\bra\cdot,\cdot\ket_a$ the scalar product in $L^2(a,\infty)$. Then 
$\bra f , H_m^\theta f\ket = \lim_{a\to0} \bra f,H_m^\theta f\ket_a$ and
\[
\bra f,H_m^\theta f\ket_a = \bra f,-\i(\partial_x +(1/2-m)Q^{-1})A_{1/2-m}^\max f\ket_a = \bra A_{1/2-m}^\max f , A_{1/2-m}^\max f\ket_a +\i\bar{f}(a) A_{1/2-m}^\max f(a).
\]
On a neighborhood of the origin we have
\begin{eqnarray*}
\i A_{1/2-m}^\max f(x) & = & \left(\partial_x +\frac{m-1/2}{x}\right)
\left(c x^{1/2-m}\cos\theta +c x^{1/2+m}\sin\theta +f_0(x) \right)\\
 & = & \left(\partial_x +\frac{m-1/2}{x}\right) \left(c x^{1/2+m}\sin\theta +f_0(x) \right)
= 2mc\sin\theta  x^{m-1/2} +o(\sqrt{x})
\end{eqnarray*}
by Proposition \ref{pr:lm0}. Then by the same proposition we get
\begin{eqnarray*}
\i\bar{f}(x) A_{1/2-m}^\max f(x) & = & (\bar{f}_0(x)+\bar{c}u_\theta(x)) 
(2mc\sin\theta  x^{m-1/2} +o(\sqrt{x})) \\
 & = & 2m|c|^2\sin\theta \left(x^{1/2-m}\cos\theta + x^{1/2+m}\sin\theta \right)
x^{m-1/2} + o(\sqrt{x})\\
 & = & 2m|c|^2\sin\theta \cos\theta +o(1).
\end{eqnarray*}
Hence $\ds \lim_{a\to0}\i\bar{f}(a) A_{1/2-m}^\max f(a)=m|c|^2\sin2\theta$.
\qed

\begin{proposition}\label{pr:fk}
Let $0<m<1$. Then $L_m^\m$ is a positive hermitian operator with deficiency indices $(1,1)$. The operators $H_m=H_m^{\pi/2}$ and $H_{-m}=H_m^{0}$ are respectively the Friedrichs and the Krein
extensions of $L_m^\m$. If $0\leq\theta\leq\pi/2$ then $H_m^\theta$ is a positive self-adjoint extension of $L_m^\m$. If $\pi/2<\theta<\pi$ then the self-adjoint extension $H_m^\theta$ of $L_m^\m$ has exactly one strictly negative eigenvalue and this eigenvalue is simple.
\end{proposition}

\proof We have, by Hardy inequality and Proposition \ref{pr:lm0}, $L_m^\m\geq m^2Q^{-2}$ as quadratic forms on $H_0^1$, so $L_m^\m$ is positive. The operators $H_m^\theta$ have the same form domain if
$\theta\neq\pi/2$, namely $\cD(A_{1/2-m}^\max)$, and $H_m^{\pi/2}$ has $H_0^1$ as form domain, which is strictly smaller.

Thus to finish the proof it suffices to show the last assertion of the proposition. Recall the modified Bessel function $K_m$ (see (\ref{eq:modifbessel})). It is easy to see that
$u_{m,k}:=\sqrt{kx}K_m(kx)$ solves $L_m^\max u_{m,k}=k^2u_{m,k}$.
Using (\ref{eq:besselfunc}), one gets that
\[
u_{m,k}\sim \frac{\pi}{2\sin\pi  m}\left(\frac{1}{\Gamma(1-m)}(kx/2)^{-m+1/2}-
\frac{1}{\Gamma(1+m)}(kx/2)^{m+1/2}\right)
\]
so that if $(k/2)^2m=-\tan\theta\Gamma(1+m)/\Gamma(1-m)$, then $u_{m,k}\in\cD(L_m^\theta)$. This proves that $L_m^\theta$ has a negative eigenvalue for $\pi/2<\theta<\pi$. It cannot have more eigenvalues, since $L_m^\min$ is positive and its deficiency indices are just $(1,1)$.
\qed

\begin{remark}{\rm The fact that $H_{\pm m}$ are the Friedrichs and the Krein extensions of $L_m^\m$ also follows from Proposition \ref{pr:KF} because we know that these are the only homogeneous extensions of $L_m^\m$.}
\end{remark}


\begin{proposition}\label{pr:f0k}
$L_0^\m$ is a positive hermitian operator with deficiency indices $(1,1)$ and its Friedrichs and Krein extensions coincide and are equal to $H_0=H_0^{\pi/2}$. The domain of $H_0$ is a dense subspace of $\cD(A_{1/2})$ and for $f\in\cD(H_0)$ we have $\bra f,H_0 f\ket=\|A_{1/2}f\|^2$. Thus the quadratic form of $H_0$ equals $A_{1/2}^*A_{1/2}$. If $0\leq\theta <\pi$ and $\theta \neq \pi/2$ then the self-adjoint extension $H_0^\theta$ of $L_0^\m$ has exactly one strictly negative eigenvalue.
\end{proposition}

\proof Since $L_0^\m$ has only one homogeneous self-adjoint extension, this follows from Proposition \ref{pr:KF} and Remark \ref{re:KF}. For the assertions concerning the quadratic form it suffices to apply
Proposition \ref{pr:mM}.
\qed

We can summarize our results in the following theorem:
\begin{theorem}\label{thm:hmfact} Let $m>-1$. Then the operators $H_m$ are positive, self-adjoint,
homogeneous of degree $2$ with $\sp H_m=\bar\R_+$. Besides we have the
following table:

\begin{tabular}{llc}
$m\geq 1$: & $H_m=A_{1/2+m}^*A_{1/2+m} \ = \ A_{1/2-m}^*A_{1/2-m},$ & $H_0^1=\cQ(H_m)$, \\[1ex]
 & & $H_m=L_m^\min=L_m^\max$; \\[2ex]
$0<m<1$: & $H_m=A_{1/2+m}^*A_{1/2+m} \ = \ \left(A_{1/2-m}^{\min}\right)^*A_{1/2-m}^{\min}$ & $H_0^1=\cQ(H_m)$, \\[1ex]
 & & \hskip -5ex $H_m$ is the Friedrichs ext. of $L_m^\min$; \\[2ex]
$m=0$: & $H_0=A_{1/2}^*A_{1/2}$, & $H_0^1 +\C\xi\zeta_0$ dense in $\cQ(H_0)$, \\[1ex]
 & &  \hskip -8ex $H_0$ is the Friedrichs and Krein ext. of $L_0^\min$; \\[2ex]
$-1<m<0$: & $H_m=\left(A_{1/2+m}^{\max}\right)^*A_{1/2+m}^{\max}$, & $H_0^1                         +\C\xi\zeta_m=\cQ(H_m)$, \\[1ex]
 & &  $H_m$ is the Krein ext. of $L_m^\min$.
\end{tabular}
\end{theorem}

In the region $-1<m<1$ (which is the most interesting one), it is quite remarkable that for strictly positive $m$ one can factorize $H_m$ in two different ways, whereas for $m\leq0$ only one
factorization appears.

As an example, let us consider the case of the Laplacian $-\partial_x^2$, i.e. $m^2=1/4$. The operators $H_{1/2}$ and $H_{-1/2}$ coincide with the Dirichlet and Neumann Laplacian respectively. One usually factorizes them as $H_{1/2}=P^*_{\min}P_{\min}$ and $H_{-1/2}=P^*_{\max}P_{\max}$, where $P_{\min}$ and $P_{\max}$ denote the usual momentum operator on the half-line with domain $H_0^1$ and $H^1$ respectively. The above analysis says that, whereas for the Neumann Laplacian this is the only factorization of the form $S^*S$ with $S$ homogeneous, in the case of the Dirichlet Laplacian one can also factorize it in the rather unusual following way
\[
H_{1/2}=\left(P_{\min}+\i Q^{-1}\right)^* \left(P_{\min}+\i Q^{-1}\right).
\]

\bep\label{re:mon}
The family $H_m$ has the following property:
\begin{eqnarray*}
0\leq m \leq m'& \Rightarrow & H_{m}\leq H_{m'},\\
0\leq m < 1 & \Rightarrow & H_{-m}\leq H_{m}.
\end{eqnarray*}
\eep


\subsection{The non hermitian case: numerical range and dissipativeness}

In this section we come back to the non hermitian case. We study the numerical
range of the operators $H_m$ in terms of the parameter $m$. As a consequence
we obtain dissipative properties of $H_m$. 

\begin{proposition} Let $m\neq 0$.
\begin{compactenum}
 \item[i)] If $0\leq\arg m\leq\pi/2$, then
   $\Num(H_m)=\{z\ |\ 0\leq\arg  z\leq 2\arg m\}$. Hence $H_m$ is
   maximal sectorial and $\i H_m$ is dissipative. 
 \item[ii)] If $-\pi/2\leq\arg m\leq 0$, then
   $\Num(H_m)=\{z\ |\ 2\arg m \leq \arg z \leq 0 \}$. Hence $H_m$ is
   maximal sectorial and $-\i H_m$ is dissipative. 
 \item[iii)] If $|\arg m|\leq \pi/4$, then $-H_m$ is dissipative. 
 \item[iv)] If $\pi/2 < |\arg m|<\pi$, then $\Num(H_m)=\C$.
\end{compactenum}
\end{proposition}

\begin{remark}{\rm For $m=0$ and $\arg m=\pi$, $H_m$ is selfadjoint so
    that $\Num(H_m)=\sp(H_m)=[0,+\infty[$. }
\end{remark}

\proof First note that since $H_m$ is homogeneous, if a point $z$ is
in the numerical range $\R_+z$ is included in the numerical
range. Thus the numerical range is a closed convex cone. Moreover,
since $H_m^*=H_{\bar{m}}$ it suffices to consider the case $\Im(m)>
0$.

Let us recall that for $\Re m>-1$ the operator $H_m$ is defined by:
\[
H_m f = -f'' +(m^2-1/4)x^{-2}f, \quad
f\in\cD(H_m)=\cD(L_m^\m)+\C\xi\zeta_m.  
\] 
Thus $C_\c^\infty +\C\xi\zeta_m$ is a core for $H_m$. Let $0<a<1$,
$c\in\C$, and $f$ a function of class $C^2$ on $\R_+$ such that
$f(x)=cx^{m+1/2}$ for $x<a$ and $f(x)=0$ for large $x$. By what we
just said the set of functions of this form is a core for $H_m$. We
set $V(x)=(m^2-1/4)x^{-2}$ and note that for any $f\in\cD(H_m)$:
\begin{align*}
\bra f, H_m f \ket & =\lim_{b\to0}\int_b^\infty
\big(-(\bar{f}f')'+|f'|^2 +V|f|^2\big) \d x \\ 
 & =\lim_{b\to0}\left(\bar{f}(b)f'(b)+ \int_b^\infty\big(|f'|^2
+V|f|^2\big) \d x \right). 
\end{align*}
If $f$ is of the form indicated above we have
$\bar{f}(b)=\bar{c}b^{\bar{m}+1/2}$ and $f'(b)=(m+1/2)cb^{m-1/2}$
for $b<a$ hence $\bar{f}(b)f'(b)=|c|^2(m+1/2)b^{2\Re m}$. To
simplify notations we set $m=\mu+\i\nu$ with $\mu,\nu$ real. Thus we
get
\begin{align*}
\bra f, H_m f \ket & = \lim_{b\to0}\left(|c|^2(m+1/2)b^{2\mu}+
\int_b^\infty\big(|f'|^2 +V|f|^2\big) \d x \right) \\ 
 & =\lim_{b\to0}\left(|c|^2(m+1/2)b^{2\mu}+ \int_b^a\big(|f'|^2
+V|f|^2\big) \d x \right) +\int_a^\infty\big(|f'|^2 +V|f|^2\big) \d x.
\end{align*}
But for $b<a$ we have
\begin{align*}
\int_b^a\big(|f'|^2 +V|f|^2\big) \d x & =
|c|^2\int_b^a\big(|m+1/2|^2 x^{2\mu-1}+ 
(m^2-1/4)x^{-2}x^{2\mu+1}\big) \d x \\
 & = |c|^2 (m+1/2)\int_b^a ( \bar{m}+1/2+ m-1/2)  x^{2\mu-1} \d x \\
 & = |c|^2 (m+1/2)\int_b^a (x^{2\mu})' \d x \ =\  |c|^2
(m+1/2)\big( a^{2\mu} - b^{2\mu} \big). 
\end{align*}
Thus we get
\begin{equation}\label{eq:numrange}
\bra f, H_m f \ket = |c|^2 (m+1/2) a^{2\mu} +
\int_a^\infty\big(|f'|^2 +V|f|^2\big) \d x =: \Psi(a,c,f).
\end{equation}
So the numerical range of $H_m$ coincides with the closure of the
set of numbers of the form $\Psi(a,c,f)$ with $0<a<1$, $c\in\C$, and
$f$ a function of class $C^2$ on $x\geq a$ which vanishes for large
$x$ and such that the derivatives $f^{(i)}(a)$ coincide with the
corresponding derivatives of $cx^{m+1/2}$ at $x=a$ for $0\leq i \leq
2$. The map $f\mapsto \int_a^\infty\big(|f'|^2 +V|f|^2\big) \d x$ is
continuous on $H^1(]a,+\infty[)$, the functions
of class $C^2$ on $[a,\infty[$ vanishing for large $x$ are dense
in this space, and the functionals $f\mapsto f'(a)$ and
$f\mapsto f''(a)$ are not continuous in the $H^1$ topology.
Hence  we may consider in the definition of $\Psi(a,c,f)$
functions $f\in H^1(]a,+\infty[)$ such that $f(a)=c a^{m+1/2}$
without extending the numerical range. 

Let $\gamma <\frac{1}{2}$, $\delta<-\frac{1}{2}$ and $R>a$, and let
$$
f(x)=\left\{\begin{array}{lll} x^{m+1/2} & {\rm if} &
x<a,\\ a^{m+1/2-\gamma}x^\gamma & {\rm if} & a\leq
x<R,\\ a^{m+1/2-\gamma}R^{\gamma-\delta}x^\delta & {\rm if} & R\leq
x.\end{array} \right.  
$$
Then one may explicitly compute
\begin{eqnarray*}
\lefteqn{(m+1/2) a^{2\mu}+ 
\int_a^\infty\big(|f'|^2 +V|f|^2\big) \d x}\\  
& = & \frac{a^{2\mu}}{1-2\gamma}(m+1/2-\gamma)^2
+a^{2\mu+1-2\gamma}R^{2\gamma-1}
\left(\frac{\delta^2+m^2-1/4}{1-2\delta}-
\frac{\gamma^2+m^2-1/4}{1-2\gamma} 
\right). 
\end{eqnarray*}
For $\gamma<\frac12$, the argument of the first term is
$2\arg(m+\frac12-\gamma)$ and the second term vanishes as
$R\to+\infty$. Using the fact that the numerical range is a convex
cone, we thus have: 
\begin{enumerate}
 \item If $\mu\geq 0$, then 
$\{z \ |\ 0\leq\arg  z\leq 2\arg m \}\subset \Num(H_m)$,  
\item If $-1<\mu<0$ then $\Num(H_m)=\C$.
\end{enumerate}
It remains to prove the reverse inclusion of 1.

We first consider the case $\mu>0$. Observe that in
\eqref{eq:numrange} $a$ may be taken as small as we wish. Hence we
can make $a\to0$ and we get  
\[
\bra f, H_m f \ket =\int_0^\infty\big(|f'|^2 +V|f|^2\big) \d x
=\|Pf\|^2+(m^2-1/4)\|Q^{-1}f\|^2,
\]
and the result follows from Proposition \ref{prop:hardy}.

On the other hand, if $\mu=0$ then the formula is different:
\[
\bra f, H_m f \ket= (m+1/2)|c(f)|^2+\|Pf\|^2+(m^2-1/4)\|Q^{-1}f\|^2,
\]
where $c(f)=\lim_{x\to0} x^{-(m+1/2)} f(x)$ is a continuous linear
functional on $\cD(H_m)$ which is nontrivial except in the case
$m=0$, cf. \eqref{eq:soL} and \eqref{eq:soL0}. In particular we have
\[
\Im \bra f, H_m f \ket = \nu \left(|c|^2 a^{2\mu}
+2\mu\int_a^\infty x^{-2}|f|^2 \d x \right)\geq 0.
\]
Since we have established the last two identities for $f$ in a core
of $H_m$,  they remain valid on $\cD(H_m)$.  \qed

As a last result, let us mention that the factorization obtained in Theorem \ref{thm:hmfact} can be extended to the complex case (see also \eqref{eq:hmfactor}), and can thus be
used as an alternative definition of $H_m$:

\begin{proposition}\label{prop:domhm} For $\Re m >-1$ we have
\begin{eqnarray*}
\cD(H_m)&:=&\left\{f\in\cD(A_{m+\frac12}^{\max})\ |\  
A_{m+\frac12}^{\max}f\in\cD(A_{\overline{m}+\frac12}^{\max*})\right\},\\
H_m f&:=& A_{\overline{m}+\frac12}^{\max*}A_{m+\frac12}^{\max} f,
\hspace{2mm} f\in\cD(H_m).
\end{eqnarray*} 
\end{proposition}

\proof Using Proposition \ref{pr:mM} and \ref{re:domt} we have $\cD(H_m)\subset \left\{ f\in\cD(A_{m+\frac12}^{\max})\ |\ A_{m+\frac12}^{\max}f\in\cD(A_{\overline{m}+\frac12}^{\max*})\right\}$.
One then prove the reverse inclusion using Proposition \ref{pr:mM} and \ref{df:hm}.
\qed

\section{Spectral projections of $H_m$ and the Hankel transformation} \label{s:hank}

In this section, we provide an explicit spectral representation of the operator $H_m$ in terms of Bessel functions.

Recall that the (unmodified) Bessel equation reads
$$
x^2w''(x)+xw'(x)+(x^2-m^2)w=0.
$$
It is well known that the Bessel function of the first kind, $J_m$ and $J_{-m}$ (see \eqref{eq:bessel}), solve this equation. Other solutions of the Bessel equations are the so-called Bessel functions of the third kind (\cite{W}) or the Hankel functions:
$$
H^\pm_m(z)=\frac{J_{-m}(z)-\e^{\mp\i m\pi}J_m(z)}{\pm\i\sin(m\pi)}.
$$
(When $m$ is an integer, one replaces the above expression by their limits). 
We have the relations
\[
J_m(x)=\e^{\pm\i\pi\frac{m}{2}}I_m(\mp\i x),\ \ \ 
H^\pm(x)=\mp\frac{2\i}{\pi}\e^{\mp\i\pi\frac{m}{2}}K_m(\mp\i x).
\]


We know that $H_m$ has no point spectrum. Hence, for any $a<b$ the Stone formula says
\begin{equation}\label{eq:stone}
\one_{[a,b]}(H_m)=\slim_{\epsilon\searrow0} \frac{1}{2\pi\i}\int_a^b
\left(G_m(\lambda+\i\epsilon)-G_m(\lambda-\i\epsilon)\right)\d\lambda.
\end{equation}

Using (\ref{eq:reso}) we can express the boundary values of the
integral kernel of the resolvent at $\lambda\in\, ]0,\infty[$ by
solutions of the standard Bessel equation:
\begin{equation*}\label{reso1}
G_m(\lambda\pm\i0;x,y) \ := \ \lim_{\epsilon\searrow0}
G_m(\lambda\pm\i\epsilon;x,y) \ = \ \left\{\begin{array}{lcc}
\pm\frac{\pi\i}{2}\sqrt{xy}J_m(\sqrt\lambda x)H_m^\pm(\sqrt \lambda
y) & {\rm if} & x<y, \\ \pm\frac{\pi\i}{2} \sqrt{xy}J_m(\sqrt\lambda
y)H_m^\pm(\sqrt\lambda x) & {\rm if} & 
x>y.  \end{array}\right.
\end{equation*}
Now 
\begin{eqnarray*}
 & & \frac{1}{2\pi\i}\left(G_m(\lambda+\i0;x,y)-G_m(\lambda-\i0;x,y)\right)\\
 & = & \left\{\begin{array}{lcc}\frac{1}{4}\sqrt{xy}J_m(\sqrt\lambda x)\left(H_m^+(\sqrt \lambda y)+
H_m^-(\sqrt\lambda y)\right)  & {\rm if} & x<y, \\
\frac{1}{4} \sqrt{xy}J_m(\sqrt\lambda y)\left(H_m^+(\sqrt \lambda
y)+H_m^-(\sqrt\lambda y)\right)  & {\rm if} & x>y;  \end{array}\right.\\
 & = & \frac{1}{2}J_m(\sqrt\lambda x) J_m(\sqrt\lambda y).
\end{eqnarray*}
Togoether with \eqref{eq:stone}, this gives an expression for the integral kernel of the spectral
projection of $H_m$, valid, say, as a quadratic form on $C_{\rm c}^\infty(\R)$.

\begin{proposition}\label{prop:hmproj} For $0<a<b<\infty$, the integral kernel of $\one_{[a,b]}(H_m)$ is
\begin{eqnarray*}
\one_{[a,b]}(H_m)(x,y) & = & \int_a^b \frac{1}{2}\sqrt{xy}J_m(\sqrt\lambda x) J_m(\sqrt\lambda y)\d\lambda\label{hankel1}\\
 & = & \int_{\sqrt a}^{\sqrt b} \sqrt{xy}J_m(k x) J_m(k y)k\d k.\label{hankel2}
\end{eqnarray*}
\end{proposition}

Let $\cF_m$ be the operator on $L^2(0,\infty)$ given by
\beq\label{def:fm}
\cF_m:f(x)\mapsto \int_0^\infty J_m(kx)\sqrt{kx}f(x)\d x
\eeq
Up to an inessential factor, $\cF_m$ is the so-called Hankel transformation.

\bet \label{th:hank} $\cF_m$ is a unitary involution on
$L^2(0,\infty)$ diagonalizing $H_m$, more precisely
$$
\cF_mH_m\cF_m^{-1}=Q^2.
$$
It satisfies $\cF_m\e^{\i tD}=\e^{-\i tD}\cF_m$ for all $t\in\R$.
\eet

\proof Obviously, $\cF_m$ is hermitian. Proposition \ref{prop:hmproj} can be rewritten as
\[
\one_{[a,b]}(H_m)=\cF_m\one_{[a,b]}(Q^2)\cF_m^*.
\]
Letting $a\to0$ and $b\to\infty$ we obtain $\one=\cF_m\cF_m^*$.  This implies that
$\cF_m$ is isometric. Using again the fact that it is hermitian we see that it is unitary.
\qed

\section{Scattering theory of $H_m$} \label{s:scm}

For the sake of completeness we give a short and self-contained
description of the scattering theory for the operators $H_m$ with
real $m$. 

\bet If $m,k>-1$ are real then the wave operators associated to the
pair $H_m,H_k$ exist and
\[
\Omega_{m,k}^\pm :=\lim_{t\to\pm\infty}\e^{\i tH_m}\e^{-\i tH_k}
=\e^{\pm \i(m-k)\pi/2}\cF_m\cF_k.
\] 
In particular the scattering operator $S_{m,k}$ for the pair
$(H_m,H_k)$ is a scalar operator
$S_{m,k}=\e^{\i\pi(m-k)}\one$.  
\eet 
\proof
Note that $\Omega^\pm_{m,k}:=\e^{\pm \i(m-k)\pi/2} \cF_m \cF_k$ is a
unitary operator in $L^2$ such that $\e^{-\i
  tH_m}\Omega^\pm_{m,k}=\Omega^\pm_{m,k}\e^{-\i tH_k}$ for all $t$.
Thus to prove the theorem it suffices to show that $(\Omega^\pm_{m.k}
- 1) \e^{-\i tH_k}\to 0$ strongly as $t\to\pm\infty$. Let $\pi_a$ be
the operator of multiplication by the characteristic function of the
interval $]0,a[$ and $\pi_a^\perp=1-\pi_a$.  Then from Theorem
\ref{th:hank} it follows easily that $\pi_a\e^{-\i tH_m}\to0$ and
$\pi_a\e^{-\i tH_k}\to0$ strongly as $t\to\pm\infty$ for any
$a>0$. Thus we are reduced to proving
\[
{\textstyle\lim_{a\to\infty}\sup_{\pm t>0}} \|\pi_a^\perp
(\Omega^\pm_{m,k} - 1) \e^{-\i tH_k}f \| = 0 
\quad \text{ for all } f\in L^2. 
\]
By using again Theorem \ref{th:hank} we get
\[
(\Omega^\pm_{m,k} - 1) \e^{-\i tH_k}=
\e^{\mp\i k\pi/2} ( \e^{\pm\i m\pi/2}\cF_m -\e^{\pm\i k\pi/2}\cF_k )
\e^{-\i t Q^2 } \cF_k
\]
hence it will be sufficient to show that
\begin{equation}\label{eq:app}
{\textstyle\lim_{a\to\infty}\sup_{\pm t>0}}
\|\pi_a^\perp( \e^{\pm\i k\pi/2}\cF_k -\e^{\pm\i m\pi/2}\cF_m )
\e^{-\i t Q^2 }g\|=0
\quad \text{ for all } g\in C_\c^\infty(\R_+). 
\end{equation}
Let us set $j_m(x)=\sqrt{x}J_m(x)$ and $\tau_m =m\pi/2+\pi/4$. Then 
$(\cF_m h)(x)=\int_0^\infty j_m(xp) h(p) \d p$ and from the
asymptotics of the Bessel functions we get
\begin{equation}\label{def:jinvol}
{\textstyle\sqrt{\frac{\pi}{2}}}j_m(y)=
\cos(y-\tau_m) +j_m^\circ(y) \quad \text{ where }
j_m^\circ(y) \sim O(y^{-1}). 
\end{equation}
If we set $g_t(p)=(\pi/2)^{1/2}\e^{-\i tp^2}g(p)$  and
$G_t^\pm= ( \e^{\pm\i k\pi/2}\cF_k -\e^{\pm\i m\pi/2}\cF_m )g_t$
\[
G_t^\pm(x)  =\int 
(\e^{\pm\i k\pi/2}\cos(xp-\tau_k) - 
\e^{\pm\i m\pi/2}\cos(xp-\tau_m)) g_t(p) \d p 
+ \int (j^\circ_k(xp)- j^\circ_m(xp)) g_t(p) \d p.
\]
The second contribution to this expression is obviously bounded by
a constant time $ |x|^{-1} \int |g_t(p)/p| \d p$ and the $L^2(\d x)$
norm of this quantity over $[a,\infty[$ is less than $Ca^{-1/2}$ for
some number $C$ independent of $t$. Thus we may forget this term in
the proof of \eqref{eq:app}. 

Finally, we consider the first contribution to $G_t^+$ for example:
since 
\[
\e^{\i k\pi/2}\cos(xp-\tau_k) - \e^{\i m\pi/2}\cos(xp-\tau_m)=
\e^{-\i xp +\i\pi/4}  (\e^{\i k\pi} - \e^{\i m\pi})/2
\]
we get an integral of the form $\int \e^{-\i p(xp+tp) } g(p) \d p$
which is rapidly decaying in $x$ uniformly in $t>0$ because
$g\in C_\c(\R_+)$ and there are no points of stationary phase.
This finishes the proof of \eqref{eq:app}.
\qed

Since $H_m$ and $H_k$ are homogeneous of degree $-2$ with respect to
the operator $D$, which has simple spectrum, we may apply Proposition
\ref{pr:wv} with $A=D$ and deduce that the wave operators are
functions of $D$. Our next goal is to give explicit formulas for these
functions.

Let $\cJ:L^2\to L^2$ be the unitary involution
\[ 
\cJ f(x)=\frac{1}{x}f(\frac{1}{x}).
\] 
Clearly $\cJ\e^{\i\tau D}=\e^{-\i\tau D}\cJ$ for all $\tau\in\R$ and $\cJ Q^2\cJ=Q^{-2}$. In particular, the operator
\begin{equation}\label{eq:gm}
\cG_m:=\cJ\cF_m
\end{equation}
is a unitary operator on $L^2$ which commutes with all the $\e^{\i\tau D}$. Hence there exists $\Xi_m:\R\to\C,$ $|\Xi_m(x)|=1$ a.e. and $\cG_m=\Xi_m(D)$. Moreover we have
\[
\cF_m\cF_k=\cF_m\cJ\cJ\cF_k=\cG_m^*\cG_k,
\]
so that
\[
\Omega_{m,k}^\pm=\e^{\pm \i(m-k)\pi/2}\cG_m^*\cG_k=\e^{\pm \i(m-k)\pi/2}\frac{\Xi_k(D)}{\Xi_m(D)}.
\]
Note that $\cG_mH_m\cG_m^*=\cJ Q^2\cJ=Q^{-2}$.

\begin{theorem}\label{thm:wavehm} For $\Re(m)>-1$,
\[
\cG_m=\e^{\i\ln(2)D} \frac{\Gamma(\frac{m+1+\i D}{2})}{\Gamma(\frac{m+1-\i D}{2})}.
\]
Therefore, for $\Re(m)$ and $\Re(k)>-1$, the wave operators for the pair $(H_m,H_k)$ are equal to
$$
\Omega_{m,k}^\pm = \e^{\pm \i(m-k)\pi/2} \frac{\Gamma(\frac{k+1+\i D}{2}) \Gamma(\frac{m+1-\i D}{2})}{\Gamma (\frac{k+1-\i D}{2})\Gamma(\frac{m+1+\i D}{2})}.
$$ 
\end{theorem}

For the proof we need the following representation of Bessel
functions: 

\begin{lemma}\label{lm:wat}
For any $m$ such that $\Re(m)>-1$ the following identity holds in the
sense of distributions:
$$
J_m(x)=\frac{1}{4\pi}\int_{-\infty}^{+\infty} \frac{\Gamma(\frac{m+\i
    t+1}{2})}{\Gamma(\frac{m-\i t+1}{2})}\left(\frac{x}{2}\right)^{-\i
  t-1}\d t. 
$$
\end{lemma}
\proof If $\Re(m)>0$ one has the following representation of the
Bessel function $J_m(x)$, cf. \cite[ch. VI.5]{W}: 
\begin{eqnarray}
J_m(x) & = & \frac{1}{2\pi\i}\int_{c-\i\infty}^{c+\i\infty}
\frac{\Gamma(z)}{\Gamma(m-z+1)}
\left(\frac{x}{2}\right)^{m-2z}\d z \nonumber\\ 
 & = & \frac{1}{4\pi} \int_{-\infty}^{+\infty}
\frac{\Gamma(c+\i\frac{t}{2})}{\Gamma(m+1-c-\i\frac{t}{2})}
\left(\frac{x}{2}\right)^{m-2c-\i t} \d t, \label{barnes} 
\end{eqnarray}
where $\ds c\in\, \left]0,\frac{\Re m}{2}\right[$. Note that the subintegral function is everywhere
analytic except for the poles at $z=0,-1,-2,\dots$, all of them on the left hand side of the contour. By the Stirling asymptotic formula, the subintegral function can be estimated by $|z|^{-1-\Re m+2c}$ at infinity, hence it is integrable.

We shall extend the formula (\ref{barnes}) for $\Re m>-1$ and 
$c\in ]0,\Re(m)+1[$. For that purpose we have to understand it in the
distributional sense, that is after smearing it with a function of
 $x$ belonging to $C_{\rm c}^\infty$.

Let $\varphi\in C_{\rm c}^\infty$ and $\ds \phi(z):=\frac{1}{4\pi}\int_0^{+\infty} \left(\frac{x}{2}\right)^z \varphi(x)\d x.$ For $\Re m>0$ and $0<c<\frac{\Re m}{2}$ we thus have
\begin{equation}\label{eq:barnesdistrib}
\int_0^\infty J_m(x)\varphi(x)\d x=\int_{-\infty}^{+\infty} \frac{\Gamma(c+\i\frac{t}{2})}{\Gamma(m+1-c-\i\frac{t}{2})} \phi(m-2c-\i t) \d t. 
\end{equation}

Since $\varphi\in \ccf$ the function $\phi$ is holomorphic and for any $K\subset \C$ compact and $n\in\N$ there exists $C_{K,n}$ s.t.
\begin{equation}\label{eq:barnesestim1}
|\phi(z+\i t)|\leq C_{K,n} \bra t\ket^{-n}, \quad \forall z\in K, \,
\forall t\in \R,  
\end{equation}
where $\bra t\ket =\sqrt{1+t^2}$. Likewise, the function $\ds z\mapsto
\theta(z)=\frac{\Gamma(z)}{\Gamma(m+1-z)}$ is holomorphic in the strip
$0<\Re(z)<\Re(m)+1$ and for any compact $K\subset\C$ there exists
$C_K>0$ such that 
\begin{equation}\label{eq:barnesestim2}
|\theta(z+\i t)|\leq C_k \bra t\ket^{2\Re(z)-\Re(m)-1}, \quad \forall z\in K, \, \forall t\in\R.
\end{equation}
Combining \eqref{eq:barnesestim1}-\eqref{eq:barnesestim2}, this proves
that the function 
$$
c\mapsto \int_{-\infty}^{+\infty}
\frac{\Gamma(c+\i\frac{t}{2})}{\Gamma(m+1-c-\i\frac{t}{2})}
\phi(m-2c-\i t) \d t 
$$ is holomorphic in the strip $0<\Re(c)<\Re(m)+1$. Moroever,
\eqref{eq:barnesdistrib} shows that this function is constant equal to
$\ds \int_0^\infty J_m(x)\varphi(x)\d x$ for $c\in \,
\left]0,\frac{\Re m}{2}\right[$. Hence \eqref{eq:barnesdistrib}
    extends to any $c$ such that $0<\Re(c)<\Re(m)+1$. In particular,
    if we chose $c=\frac{\Re(m)+1}{2}$ we get, for any $m$ with
    $\Re(m)>0$,
\begin{equation}\label{eq:barnesdistrib2}
\int_0^\infty J_m(x) \varphi(x)\d x=\frac{1}{4\pi}\int_0^\infty \d x \int_{-\infty}^{+\infty}\d t\, \frac{\Gamma(\frac{m+\i t+1}{2})}{\Gamma(\frac{m-\i t+1}{2})}\left(\frac{x}{2}\right)^{-\i t-1} \varphi(x).
\end{equation}
Using \eqref{eq:barnesestim1}-\eqref{eq:barnesestim2} once more, one
gets that the right-hand side of the above identity is holomorphic for
$\Re(m)>-1$. Since the Bessel function $J_m$ also depends on $m$ in an
holomorphic way, the left-hand side is holomorphic as well and hence
\eqref{eq:barnesdistrib2} extends to any $m$ such that $\Re(m)>-1$,
which ends the proof of the lemma.
\qed

The next lemma will also be needed. 

\begin{lemma}\label{mellin}
For a given distribution $\psi$, the operator $\psi(D)$ from $\ccf$ to $(\ccf)'$ has integral 
kernel
\begin{equation*}\label{eq:dilatkernel}
\psi(D)(x,y)=\frac{1}{2\pi\sqrt{xy}}\int_{-\infty}^{+\infty}\psi(t)\frac{y^{-\i t}}{x^{-\i t}}\d t.
\end{equation*}
\end{lemma}

\proof We use the Mellin transformation $\cM:L^2(0,\infty)\to L^2(\R)$. We
recall the formula for $\cM$ and $\cM^{-1}$:
\begin{eqnarray*}
(\cM f)(s)&:=&\frac{1}{\sqrt{2\pi}}\int_0^\infty\d x\ x^{-\frac12-\i s} f(x)\\
(\cM^{-1}g)(x)&:=&\frac{1}{\sqrt{2\pi}}\int_{-\infty}^{\infty}\d s\
x^{-\frac12+\i s}g(s).\end{eqnarray*}
The Mellin transformation diagonalizes the operator of dilations, so
that $\cM\psi(D)\cM^{-1}$ is the operator of multiplication by
$\psi(s)$.
\qed

{\noindent{\bf Proof of Theorem \ref{thm:wavehm}.}\ \ } Using
\eqref{def:fm}, \eqref{eq:gm} and Lemma \ref{lm:wat} we get that the
operator $\cG_m$ has the integral kernel
\begin{eqnarray*}
\cG_m(x,y) & = & \frac{1}{x}J_m\left(\frac{y}{x}\right)\sqrt{\frac{y}{x}}\\
 & = & \frac{1}{2\pi\sqrt{xy}}\int_{-\infty}^{+\infty}  \frac{\Gamma(\frac{m+\i t+1}{2})}{\Gamma(\frac{m-\i t+1}{2})}\left(\frac{1}{2}\right)^{-\i t}\frac{y^{-\i t}}{x^{-\i t}}\d t.
\end{eqnarray*}
Hence by Lemma \ref{mellin}, the unitary operator $\cG_m$ coincide with $\Xi_m(D)$ on $\ccf$ where
\[
\Xi_m(t)=\e^{\i\ln(2)t}\frac{\Gamma(\frac{m+1+\i t}{2})}{\Gamma(\frac{m+1-\i t}{2})}.
\]
Since $|\Xi_m(t)|=1$ for $m\in\R$, the operator $\Xi_m(D)$ is a unitary operator on $L^2$ which coincide with $\cG_m$ on the dense subspace $\ccf$, and hence $\cG_m=\Xi_m(D)$ on $L^2$.
\qed

\begin{remark}{\rm It is interesting to note that $\Xi_m(D)$ is a unitary operator for all real values of $m$ and
\beq \Xi_m^{-1}(D) Q^{-2}\Xi_m(D)\label{ham} \eeq is a function with
values in self-adjoint operators for all real $m$. $\Xi_m(D)$ is
bounded and invertible also for all $m$ such that $\Re
m\neq-1,-2,\dots$. Therefore, the formula (\ref{ham}) defines an
operator for all $\{m\mid \Re m\neq-1,-2,\dots\}\cup\R$. Clearly, for
$\Re m>-1$, this operator function coincides with the operator $H_m$
studied in this paper. Its spectrum is always equal to $[0,\infty[$
    and it is analytic in the interior of its domain.

One can then pose the following question: does this operator function
extend to a holomorphic function of closed operators (in the sense of
the definition of Subsec. \ref{ss:hol}) on the whole complex plane?}
\label{open-problem}
\end{remark}

\appendix
\section{Second order differential operators}
\label{s:so}
\setcounter{equation}{0}
\setcounter{resultcounter}{0}
\renewcommand{\theequation}{A.\arabic{equation}}
\renewcommand{\theresultcounter}{A.\arabic{resultcounter}}

To make this paper self-contained we summarize in this appendix some
facts on second order differential operators.  We are especially
interested in the case when the potential is complex and/or singular
at the origin.

\subsection{Green functions}

We consider an arbitrary complex potential $V\in L^2_{\mathrm{loc}}$ 
and a complex number $\lambda$. Let $\tilde L$ be the
distribution valued operator defined on $L^2_{\mathrm{loc}}$ by 
\beq \tilde L=-\partial_x^2+V(x).\label{ell}\eeq 
We recall that the Wronskian of two functions $f,g$ of class $C^1$ on
$\R_+$ is the function $W(f,g)$ whose value at a point $x>0$ is given
by $W_x(f,g)=f(x)g'(x)-f'(x)g(x)$. If $f,g$ are solutions of an
equation $u''=Vu$ then $W(f,g)$ is a constant which is not zero if and
only if $f,g$ are linearly independent.

We recall a standard method for constructing the Green function of a
differential operator.  An elementary computation gives:

\begin{proposition}
\label{lm:green} Suppose that $u_0$ and $u_\infty$ are
solutions of $\tilde 
Lu=\lambda u$ which are square integrable near $0$ and
$\infty$ respectively and such that $W(u_\infty,u_0)=1$. Let $g\in
L^2$ and define
\[
f_0=u_0 g_\infty+u_\infty g_0 \quad\textstyle{with}\quad
 g_0(x)=\int_0^x u_0(y)g(y)\d y, \quad g_\infty(x)=\int_x^\infty
u_\infty(y)g(y)\d y.
\]
Then the function $f_0$ satisfies $(\tilde L-\lambda)f_0=g$ and
$f_0'=u_0'g_\infty-u_\infty'g_0$. The general solution of the
equation $(\tilde L-\lambda)f=g$ can be written as $f=c_0u_0+c_\infty
u_\infty+f_0$ with $c_0,c_\infty\in \C$. We have
\begin{equation*}\label{eq:green}
f_0(x)=\int_0^\infty G(x,y)g(y)\d y  \quad \text{with} \quad
G(x,y)=\left\{\begin{array}{lll} u_0(x)u_\infty(y) & {\rm if} &
0<x<y,  \\ u_0(y)u_\infty(x) & {\rm if} & 0<y<x. \end{array}  \right.
\end{equation*}

\end{proposition}

\subsection{Maximal and minimal operators}

We denote $L_\m$ and $L_\M$ the minimal and maximal operator
associated to the differential expression (\ref{ell}). More precisely,
$L_{\M}$ is the restriction of $\tilde L$ to the space
$\cD(L_\M):=\{f\in L^2\mid \tilde L f \in L^2 \}$
considered as operator in $L^2$ and $L_\m$ is the closure of the
restriction of $L_\M$ to $C_\c^\infty$. $L_\M$ is a closed operator on $L^2$ because it is a
restriction of the continuous operator $\tilde
L:L^2_{\mathrm{loc}}\to\cD'(\R_+)$.  

>From now on  we assume that 
$\sup_{b>a}\int_b^{b+1}|V(x)| \,\d x<\infty  \text{ for }
  \text{ each } a>0$. Then we have (cf. \cite{S}):

\begin{proposition}\label{wronski}
If $f\in\cD(L_\M)$ then $f$ and $f'$ are continuous functions on
$\R_+$ which tend to zero at infinity. For $f,g\in \cD(L_{\max})$
\begin{equation}\label{eq:lag}
\lim_{x\to0} W_x(f,g) =: W_0(f,g)
\end{equation}
exists and we have 
\begin{equation}\label{eq:lagr}
\int_0^\infty (L_\max f g - f L_\max g) dx = -W_0(f,g).
\end{equation}
In particular, $W_0$ is a continuous bilinear antisymmetric form on
$\cD(L_\M)$ (equipped with the graph topology) and if one of the
functions $f$ or $g$ belongs to $\cD(L_\m)$ then $W_0(f,g)=0$.
\end{proposition}

\begin{remark}\label{re:localization}{\rm
Note that the so defined $W_0(f,g)$ depends only on the restriction
of $f$ and $g$ to an arbitrary neighborhood of zero. Hence \emph{if
  $f,g$ are continuous square integrable functions on an interval
  $]0,a[$ such that the distributions $Lf$ and $Lg$ are square
  integrable on $]0,a[$ then the limit in \eqref{eq:lagr} exists and
  defines $W_0(f,g)$.}
}\end{remark}

If $V$ is a real function the operator $L_\m$ is hermitian and
$L_\m^*=L_\M$. From \eqref{eq:lagr} we get
\begin{equation*}\label{eq:laGr}
\bra L_\M f,g \ket - \bra f, L_\M g\ket = -W_0(\bar{f},g)\equiv
\{f,g\} 
\end{equation*}
for all $f,g\in\cD(L_\M)$. Here $\{f,g\}$ is a continuous hermitian
sesquilinear form on $\cD(L_\M)$ which is zero on
$\cD(L_\m)$. Moreover, an element $f\in\cD(L_\M)$ belongs to
$\cD(L_\m)$ if and only if $\{f,g\}=0$ for all $g\in\cD(L_\M)$. A
subspace $\cE\subset\cD(L_\M)$ will be called \textit{hermitian} if it
is closed, contains $\cD(L_\m)$, and the restriction of
$\{\cdot,\cdot\}$ to it is zero. It is clear that $H$ is a closed
hermitian extension of $L_\m$ if and only if $H$ is the restriction of
$L_\M$ to a hermitian subspace.  

Now we consider the case of complex $V$.

\begin{lemma}\label{lm:dual}
Let $f\in\cD(L_\M)$. Then $f\in\cD(L_\m)$ if and only if
$W_0(f,g)=0$ for all $g\in\cD(L_\M)$. 
\end{lemma}

\proof One implication is obvious. To prove the inverse assertion let
us denote $\tilde{\bar{L}}=-\partial_x^2+\bar{V}$ acting on
continuous functions and let $\bar{L}_\m,\bar{L}_\M$ be the minimal and
maximal operators associated to $\tilde{ \bar{L}}$. It is trivial to
show that $L_\m^*=\bar{L}_\M$ hence $L_\m=\bar{L}_\M^*$ because $L_\m$
is closed. Thus $f\in L^2$ belongs to $\cD(L_\m)$ if and only if there
is $h\in L^2$ such that $\bra \bar{L}_\M g,f \ket = \bra g, h \ket$
for all $g\in\cD(\bar{L}_\M)$. But $g\in\cD(\bar{L}_\M)$ if and only
if $\bar{g}\in\cD(L_\M)$ so for $f\in\cD(L_\M)$ we get from
\eqref{eq:lagr}
\[
\bra \bar{L}_\M g,f \ket = \int_0^\infty \tilde L \bar{g} f \d x =
\int_0^\infty \bar{g} \tilde L f \d x -W_0(\bar{g},f) = \bra g,\tilde Lf \ket -W_0(\bar{g},f).
\]
Hence if $W_0(\bar{g},f)=0$ for all $g\in\cD(\bar{L}_\M)$ then $f\in\cD(L_\m)$.
\qed

We denote $\cL=\{u \mid \tilde Lu=0\}$, this is a two dimensional
subspace of $\in C^1(\R_+)$ and if $u,v\in\cL$ then $W(f,g)$ is a
constant which is not zero if and only if $u,v$ are linearly
independent. By the preceding comments, if $u\in\cL$ and
$\int_0^1|u|^2dx<\infty$ then $f\mapsto W_0(u,f)$ defines a linear
continuous form $\ell_u$ on $\cD(L_\M)$ which vanishes on
$\cD(L_\m)$. Let $L_u$ be the restriction of $L^\M$ to
$\ker\,\ell_u$. Clearly $L_u$ is a closed operator on $L^2$ such that
$L_\m\subset L_u \subset L_\M$.

\subsection{Extensions of $L_\m$}

Below by \emph{strict extension} of $L_\m$ we mean an operator $T$
such that $L_\m\subsetneq T \subsetneq L_\M$.  We denote $\xi$ a
function in $C^\infty_\c$ such that $\xi(x)=1$ for $x\leq 1$ and
$\xi(x)=0$ for $x\geq 2$.

Until the end of the subsection we assume that all the solutions of the equation $\tilde L u=0$ are square integrable at the origin.

\begin{proposition}\label{pr:sade}
$\cD(L_\m)$ is a closed subspace of codimension two of $\cD(L_\M)$ and 
\begin{equation}\label{eq:d2}
\cD(L_\m)=\{f\in\cD(L_\M) \mid W_0(u,f)=0 \ \forall u \in\cL\}=
\textstyle{\bigcap}_{u\in\cL} \ker\,\ell_u.
\end{equation}
If $u\neq0$ then $L_u$ is a strict extension of $L_\m$ and reciprocally, each strict extension of $L_\m$ is of this form. More explicitly: $\cD(L_u)=\cD(L_\m)+\C\xi u$.  We have $L_u=L_v$ if and only if $v=cu$ with $c\in\C\setminus\{0\}$.  If $V$ is real then the operator $L_\m$ is hermitian, has deficiency indices $(1,1)$, and if $u\in\cL\setminus\{0\}$ than $L_u$ is hermitian (hence self-adjoint) if and only if $u$ is real (modulo a constant factor).
\end{proposition}

\proof We first show that $\ell_u=0$ if and only if $u=0$. Indeed, if $u\neq0$ then the equation $Lv=0$ has a solution linearly independent from $u$, so that $W(u,v)\neq0$. But there is $g\in\cD(L_\M)$ such that $g=v$ on a neighborhood of zero and then $\ell_u(g)=W(u,v)\neq0$. This also proves the last assertion of the proposition.

Assume for the moment that \eqref{eq:d2} is known. If $u,v$ are linearly independent elements of $\cL$ then they are a basis of the vector space $\cL$ hence we have $\cD(L_\m)=\ker\,\ell_u\cap\ker\,\ell_v$ and so $\cD(L_\m)$ is of codimension two in $\cD(L_\M)$. Moreover, if $u\neq0$ then $\cD(L_\m)$ is of codimension one in $\ker\,\ell_u$, we have $\xi u\in\cD(L_\M)\setminus\cD(L_\m)$ and $\xi u\in\ker\,\ell_u$ hence $\cD(L_u)=\cD(L_\m)+\C\xi u$.

Then if $V$ is real the deficiency indices of $L_\m$ are $(1,1)$ because $\cD(L_\m)$ has codimension two in $\cD(L_\M)$. The space $\ker\,\ell_u$ is hermitian if and only if $\{f,f\}=0$ for all $f\in\ker\,\ell_u$. But $\ker\,\ell_u=\cD(L_\m)+\C\xi u$ so we may write $f=f_0+\lambda\xi u$ and then clearly $\{f,f\}=\{\lambda\xi u,\lambda\xi u\}=|\lambda|^2\{u,u\}=-|\lambda|^2W_0(\bar{u},u)$. So
$\ker\,\ell_u$ is hermitian if and only if $W_0(\bar{u},u)=0$. But $\bar u$ and $u$ are solutions of the same equation $Lf=0$ and $W(\bar{u},u)=W_0(\bar{u},u)=0$.  Thus $\bar u$ and $u$ must be proportional, i.e.\ there is a complex number $c$ such that $\bar u=cu$. Clearly $|c|=1$ so we may write $c=\e^{2i\theta}$ and then we see that the function $\e^{i\theta}u$ is real.

Thus it remains to prove \eqref{eq:d2} and for this we need some preliminary considerations which will be useful in another context later on. Let $v_\pm\in\cL$ such that $W(v_+,v_-)=1$. If $g$ is a function on $\R_+$ such that $\int_0^a|g|^2dx<\infty$ for all $a>0$ we set $g_\pm(x)=\int_0^xv_\pm(y)g(y)dy$.  It is easy to check that if $Lf=g$ then there is a unique pair of complex numbers $a_\pm$ such that
\begin{equation}\label{eq:sol}
f=(a_++g_-)v_++(a_--g_+)v_-
\end{equation}
and reciprocally, if $f$ is defined by \eqref{eq:sol} then $Lf=g$. Since $g_\pm'=v_\pm g$ we also have
\begin{equation}\label{eq:sold}
f'=(a_++g_-)v'_++(a_--g_+)v'_-.
\end{equation}
Now assume $h\in\cD(L_\M)$ and $W_0(u,h)=0$ for all $u\in\cL$. This is equivalent to $\ell_{v_\pm}(h)=0$. We shall prove that $W_0(f,h)=0$ for all $f\in \cD(L_\M)$ and this will imply $h\in\cD(L_\m)$ by Lemma \ref{lm:dual}. If we set $v=a_+v_+ + a_-v_-$ and $f_0=g_-v_+ - g_+v_-$ then we get $W_0(f,h)=W_0(f_0,h)$. Then
\[
W_0(f_0,h)=W_0(g_-v_+ - g_+v_-,h)=\lim_{x\to0} 
\left((g_-v_+ - g_+v_-)(x) h'(x) -(g_-v_+ - g_+v_-)'(x) h(x)\right).
\]
For a fixed $x$ we rearrange the last expression as follows:
\[
g_- v_+ h' -(g_-v_+ )' h  - g_+v_- h'  + (g_+v_-)' h = g_- W_x(v_+,h) - g_+W_x(v_-,h) - g_-'v_+h + g_+'v_- h.
\]
When $x\to 0$ the first two terms on the right hand side clearly converge to zero. The last two become $- gv_-v_+h + gv_+v_- h=0$. This finishes the proof.
\qed

\begin{remark}\label{re:reg}{\rm
If zero is a regular endpoint, i.e. $\int_0^1|V(x)|\d x<\infty$, then
for each $f\in \cD(L_\M)$ the limits $\lim_{x\to0}f(x)\equiv f(0)$ and
$\lim_{x\to0}f'(x)\equiv f'(0)$ exist. If $V$ is real we easily get
the classification of the self-adjoint realizations of $L$ in terms of
boundary conditions of the form $f(0)\sin\theta-f'(0)\cos\theta=0$.}
\end{remark}

We point out now some consequences of the preceding proof. We denote $\|h\|_x$ the $L^2$ norm of a function $h$ on the interval $]0,x[$.  Then we get $ |g_\pm(x)|\leq\|v_\pm\|_x\|g\|_x $ for all $x>0$, 
where the numbers $\|v_\pm\|_x$ are finite and tend to zero as $x\to0$. Note that in general $\|v_\pm'\|_x=\infty$ for all $x$ for at least one of the indices $\pm$. Anyway, we have
\begin{eqnarray*}
|f(x)-(a_+v_+(x)+a_-v_-(x))| & \leq &
\Big(|v_+(x)|\|v_-\|_x+|v_-(x)|\|v_+\|_x\Big)\|g\|_x \\
|f'(x)-(a_+v_+'(x)+a_-v_-'(x))| & \leq &
\Big(|v_+'(x)|\|v_-\|_x+|v_-'(x)|\|v_+\|_x\Big)\|g\|_x. 
\end{eqnarray*} 
In other terms: if $f$ is a solution of $Lf=g$ then there are complex
numbers $a_\pm$ such that as $x\to0$: 
\begin{eqnarray}
f(x) & = & a_+v_+(x)+a_-v_-(x)+
o(1)\Big(|v_+(x)|\|v_-\|_x+|v_-(x)|\|v_+\|_x\Big), \label{eq:asym} \\
f'(x) & = & a_+v_+'(x)+a_-v_-'(x)+
o(1)\Big(|v_+'(x)|\|v_-\|_x+|v_-'(x)|\|v_+\|_x\Big), 
\label{eq:asymd}
\end{eqnarray}

In the next proposition we continue to assume that all the solutions
of the equation $Lu=0$ are square integrable at the origin and keep
the notations introduced in the proof of Proposition \ref{pr:sade}.

\begin{proposition}\label{pr:domto}
A function $f\in \cD(L_\M)$ belongs to $\cD(L_\m)$ if and only if
$f=v_+g_--v_-g_+$ with $g=Lf$. In particular, if $f\in D(L_\m)$ then
for $x\to0$ we have:
\begin{equation*}\label{eq:esto}
f(x)=o(1)\Big(|v_+(x)|\|v_-\|_x+|v_-(x)|\|v_+\|_x\Big), \hspace{2mm}
f'(x)=o(1)\Big(|v_+'(x)|\|v_-\|_x+|v_-'(x)|\|v_+\|_x\Big).
\end{equation*}
\end{proposition}

\proof We take above $g=Lf$ and we get the relations \eqref{eq:sol},
\eqref{eq:sold}, \eqref{eq:asym} and \eqref{eq:asymd} for some
uniquely determined numbers $a_\pm$. If we set $v=a_+v_++a_-v_-$ and
$f_0=v_+g_--v_-g_+$ then $f=v+f_0$.  We know that $f\in D(L_\m)$ if
and only if $W_0(u,f)=0$ for all $u\in\cl$. Since $v_\pm$ form a basis
in $\cl$, it suffices in fact to have this only for $u=v_\pm$. We have
$W_0(v_\pm,f_0)=0$ because $f_0'=v_+'g_--v_-'g_+$, so that
$$
v_\pm f_0'-v_\pm'f_0=v_\pm(v_+'g_--v_-'g_+)-v_\pm'(v_+g_--v_-g_+) =
-g_\pm, 
$$ and $g_\pm(x)\to0$ as $x\to0$. Hence $W_0(v_\pm,f) =
W_0(v_\pm,v)+W_0(v_\pm,f_0) = W_0(v_\pm,v)=\pm a_\mp$, and so $f\in
D(L_\m)$ if and only if $a_\pm=0$, or if and only if $f=v_+g_--v_-g_+$
with $g=Tf$. Thus, if $f\in D(L_\m)$ then we have the relations
\eqref{eq:asym} and \eqref{eq:asymd} with $a_\pm=0$, so we have the
required asymptotic behaviours of $f$ and $f'$.  \qed

\section{Aharonov-Bohm Hamiltonian}\label{ss:ahb}
\setcounter{equation}{0}
\setcounter{resultcounter}{0}
\renewcommand{\theequation}{B.\arabic{equation}}
\renewcommand{\theresultcounter}{B.\arabic{resultcounter}}

Consider the Hilbert space $L^2(\R^2)$.  We will use simultaneously the
polar coordinates, $r,\phi$, which identify this Hilbert space with
$L^2(0,\infty)\otimes L^2(-\pi,\pi)$ by the unitary transformation
\[L^2(\R^2)\ni f\mapsto Uf\in L^2(0,\infty)\otimes L^2(-\pi,\pi)\]
given by $Uf(r,\phi)=\sqrt{r}f(r\cos\phi,r\sin\phi)$.

Let $\lambda\in\R$. We consider the magnetic hamiltonian associated to the magnetic potential $(\frac{\lambda y}{x^2+y^2},-\frac{\lambda x}{x^2+y^2})$. The curl of this potential equals zero away from the origin of coordinates and the corresponding Hamiltonian (at least for real $\lambda$) is called the Aharonov-Bohm Hamiltonian. More precisely, let $M_\lambda$ denote the minimal operator associated to the differential expression
\begin{equation}\label{def:ab}
M_\lambda:=-\left(-i\partial_x-\frac{\lambda y}{x^2+y^2}\right)^2-
\left(-i\partial_y+\frac{\lambda x}{x^2+y^2}\right)^2,
\end{equation}
a priori defined on $\ccf(\R^2\backslash \{0\})$. Clearly, $M_\lambda$ is a positive hermitian operator, homogeneous of degree $-2$. In polar coordinates, $M_\lambda$ becomes
\begin{equation}\label{eq:abpolar}
M_\lambda=-\partial_r^2+\frac{1}{r^2}\left[(-i\partial_\phi+\lambda)^2-\frac{1}{4}\right].
\end{equation}
Let $L:= -\i x\partial_y+\i y\partial_x$ be the angular
momentum. $L=-\i\partial_\phi$ in polar coordinates. Then
$L$ commutes with $M_\lambda$ (or equivalently, $M_\lambda$ is rotation
symmetric). $L$ is a self-adjoint operator with the spectrum $\sp(L)=\Z$. 
Therefore, we have a direct sum decomposition $L^2(\R^2)
=\mathop{\oplus}\limits_{l\in\Z}\cH_l$ where $\cH_l$ is the spectral
subspace of $L$ for the eigenvalue $l$. With the help of $U$ we can identify $\cH_l$ with $L^2(\R_+)$.

Using (\ref{eq:abpolar}), one immediately gets that
\begin{equation}\label{eq:abdecomp}
U M_\lambda U^* = \mathop{\oplus}\limits_{l\in\Z} L^\min_{l+\lambda}.
\end{equation}
Using general arguments, see Proposition \ref{pr:KF}, one easily gets
that the Friedrichs and the Krein extensions, denoted $M_\lambda^{FF}$
and $M_\lambda^{KK}$ respectively, of $M_\lambda$ are also homogeneous
and rotation symmetric (the reason for the double superscript will
become apparent later).

\begin{proposition}\label{prop:abdefindices}
\begin{compactenum}
\item[(i)] If $\lambda\in\Z$, then $M_\lambda$ has deficiency indices $(1,1)$. We have $M_\lambda^{FF}=M_\lambda^{KK}$, and $M_\lambda$ has no other homogeneous extension.
\item[(ii)] If $\lambda\notin \Z$, then $M_\lambda$ has deficiency indices $(2,2)$. We have $M_\lambda^{FF}\neq M_\lambda^{KK}$, and $M_\lambda$ has two other (distinct) homogeneous and rotation symmetric self-adjoint extensions $M_\lambda^{FK}$ and $M_\lambda^{KF}$.
\end{compactenum}
\end{proposition}

\begin{remark} When $\lambda\notin \Z$, $M_\lambda$ has also many homogeneous self-adjoint extensions which are not rotation symmetric.
\end{remark}

\begin{remark} If $V$ denotes the unitary operator such that
  $V=\e^{i\phi}$ in polar coordinates, then 
\begin{equation}\label{eq:absymmetry1}
V^*M_\lambda V=M_{\lambda+1}.
\end{equation} 
\end{remark}

\proof Using (\ref{eq:abdecomp}), the deficiency indices of
$M_\lambda$ are $(n,n)$ where $n=\sum_{l\in\Z} n_l$, and $(n_l,n_l)$
are the deficiency indices of $L_{l+\lambda}^\min$. By Proposition
\ref{pr:1}, we have $n_l=0$ unless $|l+\lambda|<1$ in which case
$n_l=1$. Thus, if $\lambda\in\Z$ only the term with $l=-\lambda$ has
nonzero deficiency indices, namely $n_{-\lambda}=1$, and if
$\lambda\notin \Z$ then $n_l=1$ only when $l=-[\lambda]-1$ and
$l=-[\lambda]$, where $[\lambda]$ denotes the integer part of
$\lambda$. This proves the assertions concerning the deficiency
indices.

Using (\ref{eq:absymmetry1}) we may then restrict to the case $0\leq
\lambda <1$. The result follows from the analysis of Section
\ref{ssec:lminextension}. If $\lambda=0$ the only term which is not
self-adjoint in the decomposition of $M_0$ is $L_0^\min$. Using
Proposition \ref{pr:homc} we see that $M_0$ has a unique homogeneous
self-adjoint extension. Since $M_0^{FF}$ and $M_0^{KK}$ are both
homogeneous they necessarily coincide.

We then turn to the case $0<\lambda <1$. Only the terms $L_{\lambda-1}^\min$ and $L_\lambda^\min$ are not self-adjoint. Using Proposition \ref{pr:homc} again, each of these term has exactly two homogeneous extensions $H_{\pm(\lambda-1)}$ and $H_{\pm\lambda}$ respectively, those with a $+$ sign corresponding to the Friedrichs extension and those with a $-$ sign to the Krein extension.
Hence $M_\lambda$ has 4 distinct homogeneous \emph{and} rotation symmetric self-adjoint extensions. The super-indices $FF$, $KK$, $FK$ and $KF$ correspond to the choice of the two extensions (the first index for the extension of $L_{\lambda-1}^\min$).
\qed

We can then apply the results of Section \ref{ssec:holom} to study the analiticity properties of the various homogeneous extensions of $M_\lambda$. 

\bet\label{thm:abholom} Let $n\in\Z$.
For any $\#\in\{FF,KK,FK,KF\}$ the map $]n,n+1[\ \ni\lambda\mapsto M_\lambda^\#$ extends to a holomorphic family $M^\#_z$ on the strip $\{n<\Re(z)<n+1\}.$ Moreover,
\begin{compactenum}
\item[(i)] the family $z\mapsto M_z^{FF}$ can be extended to a holomorphic family on the strip $\{ n-1<\Re(z)<n+2\}$.
\item[(ii)] the family $z\mapsto M_z^{FK}$ can be extended to a holomorphic family on the strip $\{ n-2<\Re(z)<n+1\}$.
\item[(iii)] the family $z\mapsto M_z^{KF}$ can be extended to a holomorphic family on the strip $\{ n<\Re(z)<n+3\}$.
\end{compactenum}
\eet

\proof Using Proposition \ref{prop:abdefindices}, for any $\lambda\in\ ]n,n+1[$, we have
\begin{equation}
M_\lambda^\# = \mathop{\oplus}\limits_{l\leq -n-2} H_{-l-\lambda} \oplus H_{\pm(\lambda-n-1)}\oplus H_{\pm(\lambda-n)} \mathop{\oplus}\limits_{l\geq -n+1} H_{l+\lambda}.
\end{equation}
Using Theorem \ref{th:holo}, the components $H_{-l-\lambda}$ (for
$l\leq -n-2$) have an analytic extension to the half-plane
$\Re(z)<-l+1$, the components $H_{l+\lambda}$ (for $l\geq n+1$) have
an analytic extension to the half-plane $\Re(z) > -l-1$. Similarly,
$H_{\lambda-n-1}$ (the Krein extension of $L_{\lambda-n-1}^\min$) has
an extension to the half-plane $\Re(z)>n$, $H_{-\lambda+n+1}$ to the
half-plane $\Re(z)<n+2$, $H_{\lambda-n}$ to the half-plane
$\Re(z)>n-1$ and $H_{-\lambda+n}$ to the half-plane $\Re(z)<n+1$. The
result then easily follows.  \qed

\begin{remark} The value at $z=n$ of both families $M_z^{FK}$ and
  $M_z^{FF}$   coincides with the unique homogeneous extension
  of $M_n$.
\end{remark}

\end{document}